\documentclass[a4paper,12pt,oneside]{article}

\usepackage{eucal}
\usepackage{color}
\usepackage[all]{xy}
\usepackage[centertags]{amsmath}
\usepackage{latexsym}
\usepackage{stmaryrd}
\usepackage{amsfonts}
\usepackage{amssymb}
\usepackage{amsthm}
\usepackage{fancyhdr}
\usepackage{t1enc}
\usepackage[utf8]{inputenc} 
\usepackage[english]{minitoc}
\usepackage{fancyhdr}
\usepackage{newlfont}
\usepackage[active]{srcltx}
\usepackage{graphicx}
\usepackage{t1enc}
\usepackage{tikz}
\usepackage{makeidx}
\usepackage{appendix}
\usepackage[english]{babel}			
\usepackage{float}

\theoremstyle{plain}
\newtheorem{theorem}{Theorem}[section]
\newtheorem{remark}[theorem]{Remark}
\newtheorem{lemma}[theorem]{Lemma}
\newtheorem{exercise}[theorem]{Exercise}
\newtheorem{proposition}[theorem]{Proposition}
\newtheorem{definition}[theorem]{Definition}
\newtheorem{corollary}[theorem]{Corollary}

\newtheorem{example}[theorem]{Exemple}
\newtheorem{prop}{Propri\'et\'e}[section]

\newcommand{\Z}{\mathbb Z}

\newcommand{\N}{\mathbb N}

\newcommand{\A}{\mathcal{A}}
\newcommand{\D}{\mathbb D}

\newcommand{\G}{\mathcal{G}}

\newcommand{\ds}{\displaystyle}
\def\pn {\par \noindent}
\def\md {\par \medskip}

\newcommand{\beq}{\begin{equation}}
\newcommand{\eeq}{\end{equation}}
\newcommand{\beqn}{\begin{eqnarray}}
\newcommand{\eeqn}{\end{eqnarray}}
\newcommand{\bpro}{\begin{proposition}}
\newcommand{\epro}{\end{proposition}}
\newcommand{\blem}{\begin{lemma}}
\newcommand{\elem}{\end{lemma}}
\newcommand{\bdfn}{\begin{definition}}
\newcommand{\edfn}{\end{definition}}
\newcommand{\bcor}{\begin{corollary}}
\newcommand{\ecor}{\end{corollary}}
\newcommand{\bthm}{\begin{theorem}}
\newcommand{\ethm}{\end{theorem}}
\newcommand{\bex}{\begin{example}}
\newcommand{\eex}{\end{example}}
\newcommand{\brmq}{\begin{remark}}
\newcommand{\ermq}{\end{remark}}
\newcommand{\benum}{\begin{enumerate}}
\newcommand{\eenum}{\end{enumerate}}
\newcommand{\bitem}{\begin{itemize}}
\newcommand{\eitem}{\end{itemize}}
\newcommand{\bexer}{\begin{exercise}}
\newcommand{\eexer}{\end{exercise}}
\newcommand{\bproof}{\begin{proof}}
\newcommand{\eproof}{\end{proof}}
\newcommand{\eprop} {\end{prop} }
\newcommand{\bprop}{\begin{prop}}
\theoremstyle{plain}

\usetikzlibrary{arrows}
\usetikzlibrary{decorations.markings}
\tikzstyle directed=[postaction={decorate,decoration={markings,
    mark=at position .65 with {\arrow{latex}}}}]

\providecommand{\thanks}[1]

\usepackage[english]{babel}
\title{Self-intersection on pair of pants}
\author{ElHadji Abdou Aziz Diop, Masseye Gaye}

\begin{document}
\maketitle						
\abstract

In this paper\footnote{\noindent{\bf 2010 Mathematics Subject Classification:} Primary: 32G15. Secondary 30F40.

{\bf Keywords:} closed geodesics, Schottky, self-intersection, $k$-systole.}
, we use the coding developed by R. Bowen and C. Series to compute the number of self-intersections of a closed geodesic on a pair of pants. We give lower and upper bounds on the number of self-intersections of a closed geodesic on a pair of pants. We prove a conjecture of Moira Chas and Anthony Phillips in \cite{CP}. We get also bounds for the number of closed geodesics whose self-intersection number is very close to the maximal self-intersection number on a pair of pants.

\section{Introduction}
Let $S=\D/\Gamma$  be a hyperbolic surface where $\Gamma$  is a Schottky group purely hyperbolic generated by the isometries $a_1,b_1, \cdots, a_p,b_p$. We note by $\Sigma$ the set of all bi-infinite reduced words of $\Gamma$.  The coding introduce by R.Bowen and C. Series(  see \cite{BS}, \cite{S1}, \cite{S2}) gives a bijective  correspondence between the set of closed geodesics $ \mathcal{G}^c $ of $S$ and the set $\Sigma$ under cyclic permutation; in particular a closed geodesic $\gamma\in \mathcal{G}^c$ is associated to a unique finite word  $w(\gamma)=s_1^{i_1}r_1^{j_1}\cdots s_n^{i_n}r_n^{j_n}$ under cyclic permutation where $s_i\in \{a_1,\cdots ,a_p,\bar{a}_1,\cdots ,\bar{a}_p\}$ and $r_i\in\{b_1,\cdots , b_p,\bar{b}_1,\cdots, \bar{b}_p\}$, $\bar{a}_i=a_i^{-1}$ and $\bar{b}_i=b_i^{-1}$ for all $1\leq i\leq n$.
The integer $L(\gamma)=\ds\sum_{k=1}^{n}(i_k+j_k)$ is the combinatorial length of $\gamma$.

 In this paper, we are interested in the relationship between the self-intersection and the combinatorial length of a closed geodesic on a pair of pants.
 The methods on this paper are combinatorial and we use the ideas developed by Cohen and Lustig in \cite{cohen2016paths} to compute the self-intersection number of a closed geodesic. In \cite{Bas1}, A. Basmajan study the relationship between self-intersection and the geometric length of a closed geodesic on a hyperbolic surface by using hyperbolic geometric.
 
 Let ~$\D$ denotes the Poincar\'e disk endowed with the hyperbolic metric and let $a$ and $b$ two hyperbolic isometries of $\D$ verifying the following conditions:

\benum

\begin{figure}
\begin{center}
\includegraphics[scale=0.2]{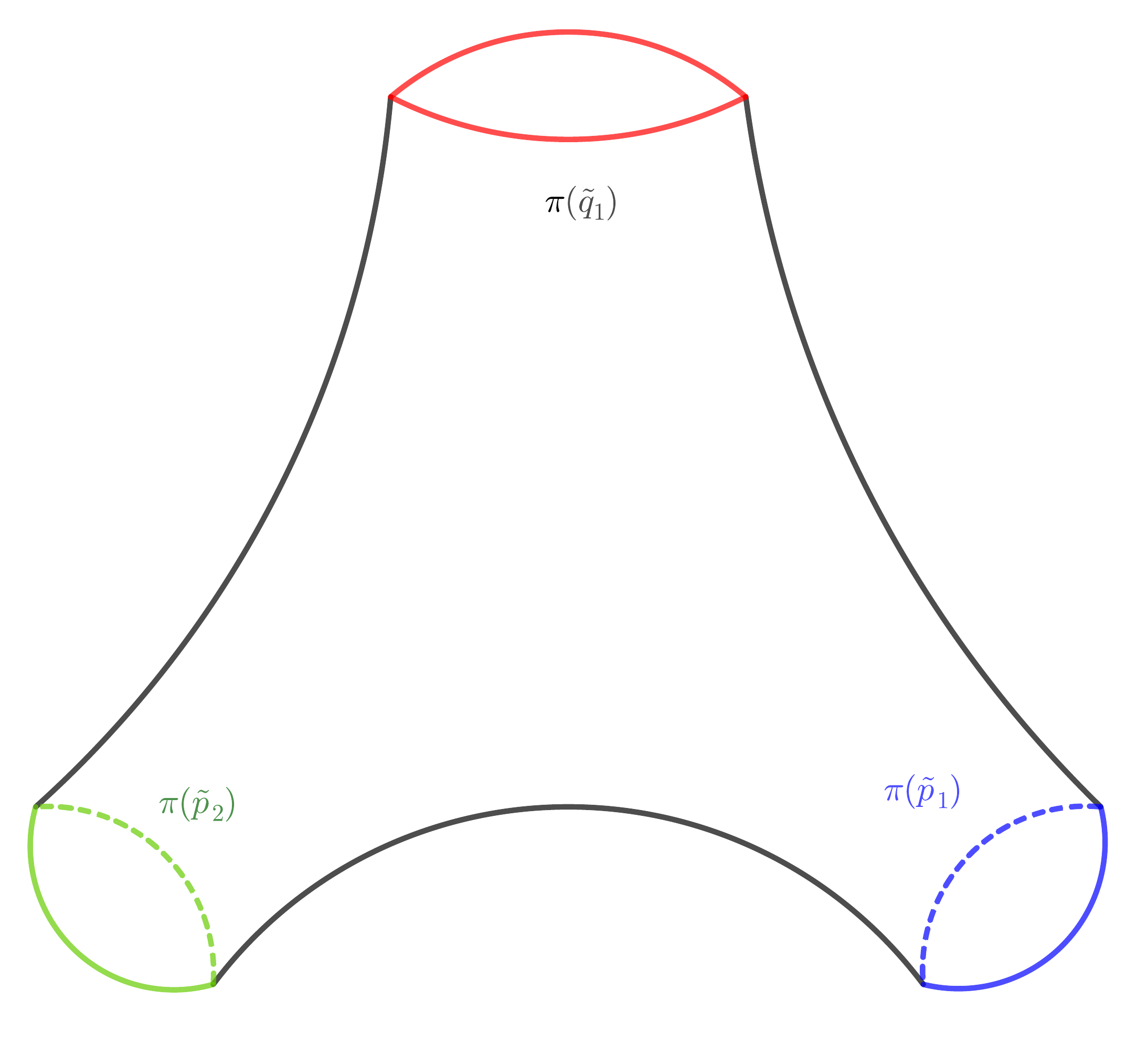}\includegraphics[scale=0.2]{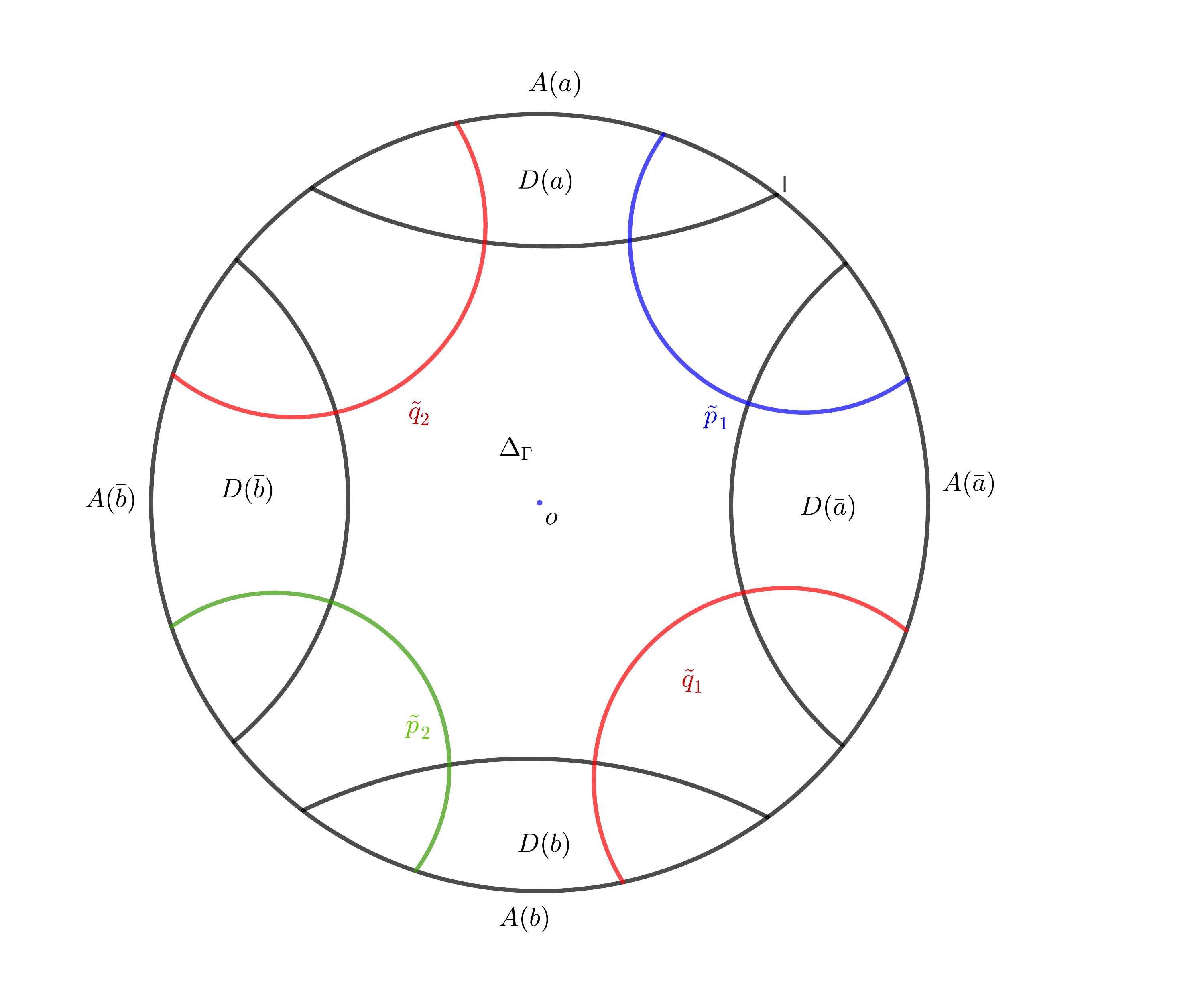}

\caption{}\label{coding}
\end{center}
\end{figure}

\item there are four disjoint euclidean closed half-disks $D(a)$, $D({\bar{a}})$, $D(b)$ and $D({\bar{b}})$ of $\D$ which are orthogonal to $\mathbb{S}^{1}$ and such that $a(D(a))=\D-D({\bar{a}})$ et $b(D(b))=\D-D({\bar{b}})$.
\item the half-disks $D(a)$ and $D({\bar{a}})$ side by side (see figure \ref{coding}).

where $\bar{a}$ and $\bar{b}$ are respectively the inverse of $a$ and $b$.
  
\eenum

Throughout this paper we will consider the configuration in figure \ref{coding}. The group $\Gamma:=\langle a, b\rangle$ generated by $a$ and $b$ is a Schottky group and the surface $P=\D/\Gamma$ is a pair of pants (see \cite{Bea}, \cite{Dal}).
Let us denote $\overline{\Gamma}=\{a,\bar{b},b, \bar{a}\}$ the set of generators of $\Gamma$ and their inverse. 
The set $\mathcal{D}=\ds\bigcap_{e\in\{a, \bar{a}, b, \bar{b}\}}\D-\overset{\circ}{D}(e)$  is a fundamental domain of $\D$ under the action of $\Gamma$.

Using this symbolic point of view, we give a method to calculate the self-intersection number of any geodesic in $ \mathcal{G}^c $.
Let $\gamma$ be a closed geodesic of $P$ and $w(\gamma)=s_1^{i_1}r_1^{j_1}\cdots s_n^{i_n}r_n^{j_n}$ the word associated to $\gamma$ where $n$, $i_1, j_1\cdots ,i_n, j_n$ are non-zero positive integers and $s_k\in \{a,\bar{a}\}$ and $r_k\in\{b,\bar{b}\}$ for all $1\leq k\leq n$.
Each cyclic permutation of $w(\gamma)$ corresponds to a unique lift $\gamma_i$  of $\gamma$ in $\D$ which crosses the fundamental domain $\Delta_{\Gamma}$ for $1\leq i\leq L(\gamma)$ . Then: $$i(\gamma;\gamma)=\#\{\gamma_i\cap\gamma_j\cap\mathcal{D}\neq\emptyset; \ \  1\leq i<j\leq L(\gamma)\}.$$
The self-intersection number $i(\gamma;\gamma)$ of $\gamma$ satisfies the following inequalities:

\bthm\label{thm1}
Let $\gamma$ be a non-simple closed geodesic of $P$ and $w(\gamma)=s_1^{i_1}r_1^{j_1}\cdots s_n^{i_n}r_n^{j_n}$ the word associated to $\gamma$, then:
$$L(\gamma)-n-1\leq i(\gamma;\gamma)\leq nL(\gamma)-n^{2}.$$
\ethm

In \cite{CP}, M.Chas and A.Phillips use the same method to tabulate self-intersection numbers for closed curves on pair of pants and they proved that the self-intersection numbers are bounded below and these bounds are sharp. They also proved that the self-intersection of a closed curve of combinatorial length $L(\gamma)$ is bounded above and this bound is sharp when $L(\gamma)$ is even and when $L$ is odd, they conjectured that the maximal self-intersection number is $\ds\frac{L^2(\gamma)-1}{4}$. The theorem \ref{thm1} enables us to prove this conjecture for closed geodesics on the pair of pants.

\bthm\label{thm2}
Let $\gamma$ be a non-simple closed geodesic of $P$ and $L(\gamma)$ combinatorial length then we have:

$$\ds\left\{\begin{array}{ll}
\ds\frac{L(\gamma)}{2}-1 \ \  \mbox{if}  \  \mbox{$L(\gamma)$ is even}\\
\\
\ds\frac{L(\gamma)-1}{2} \ \ \mbox{if} \ \mbox{ $L(\gamma)$ is odd} \\
\end{array}\right.\leq i(\gamma;\gamma)\leq  \ds\left\{\begin{array}{ll}
\ds\frac{L^2(\gamma)}{4} \ \ \mbox{if}  \  \mbox{$L(\gamma)$ is even}\\
\\
\ds\frac{L^2(\gamma)-1}{4} \ \ \mbox{if}  \ \mbox{ $L(\gamma)$ is odd} \\
\end{array}\right.$$

These bounds are sharp.

\ethm



We want to determine the proportion of closed geodesics of $P$ whose self-intersection number is very close to the maximal self-intersection number. Set 
$$\A{\ds_\epsilon}(L)=\{\gamma\in \mathcal{G}^c | L(\gamma)=L ,\ i(\gamma;\gamma)\geq (\ds\frac{1}{4}-\epsilon)L^2 \}$$
where $L$ is an integer and $\ds\epsilon$ a positive real.
 In this paper, we get the following result for a pair of pants.
\bthm\label{thm5}
 For $L$ large enough, we have:
$$\ds \frac{4\times 3^{6\ds\epsilon L}}{L}\leq\#\A{\ds_\epsilon}(L)\leq \ds\frac{2^L}{L}\big[2^2(\ds\frac{3}{2})^{\ds 6\sqrt{\epsilon}L}-3\sqrt{\epsilon}L\big]$$
and 
$$\ds\lim\limits_{L\to +\infty}\frac{\#\A_{\ds\epsilon}(L)}{\#\G^{c}(L)}=0 \ \
\mbox{for} \ \ \ds\epsilon<\ds\frac{1}{36}.$$
\ethm
The problem of counting closed geodesics has been studied in many contexts ( see \cite{Lal2}, \cite{Mar04}, \cite{Patt88}, \cite{Hub1959},\cite{Gui1986}). For more details, there is a survey of the history of this problem by Richard Sharp that was published in conjunction with Margulis's thesis in \cite{Shar04}.

Recently, there has been work on the dependence of the number of closed geodesics and their self-intersection number as well as length (see \cite{Mir} \cite{Riv12}, \cite{Ree} \cite{Sap1},\cite{Sap2} and \cite{Sap3}).

In Section 2, we recall the construction of a combinatorial model for closed geodesics on a pair of pants. We show that each geodesic can be represent as a cyclically reduced non-periodic word.
In section $3$, we use this combinatorial model to develop a method which enables us to calculate the number of self-intersections of a given closed geodesic on a pair of pants.
In section $4$, we use the formula obtained in section $3$ to prove the theorems previously enunciated.

\section{Coding closed geodesics}\label{codgeo}

Let $a$ and $b$ be two hyperbolic isometries such that we have the configuration describe in the introduction( see figure 1). 

Denote by $\Gamma=\langle a,b\rangle$ the group generated by $a$ and $b$ and by $\overline{\Gamma}=\{a,\bar{b},b, \bar{a}\}$ the set of generators of $\Gamma$ and their inverse.
 Our goal on this section is to give a way to code the closed geodesics of the hyperbolic surface $P=\D/\Gamma$. We begin to code the the points of the limit set $L_\Gamma$ of $\Gamma$.
Any geodesic of the hyperbolic surface $P$ admits a lift in the Poincaré disc which both endpoints lie in the limit set $L_\Gamma$. Together with the coding of the limit set, this implies the coding of the geodesics of $P$.

\subsection{Coding of the limit set}

In this section, we will give a coding of the points of $L_\Gamma$. For any $e\in \overline{\Gamma}$, define $A(e) $ to be the arc of $\partial \D$ inside $D(e)$ ( see figure \ref{coding} ). We will call these four arcs the first order intervals on $\partial\D$.

\bdfn

\bitem
\item A word in $\Gamma$ is a sequence $e_1e_2\cdots $ where $e_i\in \overline{\Gamma}$.

\item  A finite word $e_1e_2\cdots e_n$ is reduced if for all $1\leq i\leq n$, $e_i\neq \bar{e}_{i+1}$. It is said cyclically reduced if $e_1\neq \bar{e}_n$.

\item  A infinite or bi-infinite word is reduced if each of his finite sub-word is reduced.  
\eitem
\edfn

In order to have a better localisation of the points of the limit set and to construct the coding of the limit set we need to use the sub-arcs of the first order intervals on $\partial\D$.

\bdfn\label{interval}
For any reduced word $e_1e_2\cdots $ in $\Gamma$, the subset of $\partial\D$ defined by: 
\begin{center}
$A(e_1e_2\cdots e_m)=e_1\cdots e_{m-1}A(e_m)$
\end{center}
is called $m-th$ order interval.

\edfn

\bpro
\benum

\item For any reduced word $e_1e_2\cdots $, we have

$$A(e_1e_2\cdots e_m)\subset A(e_1e_2\cdots e_{m-1})\subset \cdots \subset A(e_1).$$

\item  If $e_1e_2\cdots e_m$ and $f_1f_2\cdots f_m$ are two reduced words such that $e_1e_2\cdots e_m\neq f_1f_2\cdots f_m$
then: $$A(e_1e_2\cdots e_m)\cap A(f_1f_2\cdots f_m)=\emptyset.$$
\eenum

\epro

\bproof
\benum
\item  We will show this point by recurrence. Because $e_1e_2\cdots $ is reduced,  $e_1\neq \bar{e}_2$ and t hen $A(e_2)\subset \partial\D\setminus A(\bar{e}_1)$. The condition on the isometries $a$ and $b$ implies that $A(e_1e_2)=e_1A(e_2)\subset A(e_1)$. We suppose that $A(e_2e_3\cdots e_m)\subset A(e_2e_3\cdots e_{m-1})$. Because $e_1$ is an isometry, calculate the image of these sets by $e_1$ gives us the result.

\item Because $e_1e_2\cdots e_m\neq f_1f_2\cdots f_m$, there exists an integer $p$ between $1$ and $m-1$ such that $e_1e_2\cdots e_p=f_1f_2\cdots f_p$ and $e_{p+1}\neq f_{p+1}$. Using the point $1$ of this proposition and the fact that $A(e_{p+1})\cap A(f_{p+1})=\emptyset$ give us the result.

\eenum
\eproof

Set $A=\ds\bigcup_{e\in\overline{\Gamma}} A(e)$ and define the following map $f:A\rightarrow \partial\D$ by $f_{|A(e)}(x)=\bar{ e}(x)$.
Let $\eta\in L_\Gamma\subset A$ be a point of the limit set then there exists $e_1\in \overline{\Gamma}$ such that $\eta\in A(e_1)$ and $f(\eta)=\bar{e}_1(\eta) \in \partial\D\setminus A(\bar{e}_1)$.
If $f(\eta) \in A$ then there exists $e_2\in \overline{\Gamma}$ such that $f(\eta)\in A(e_2)$ and in this case we have $f^2(\eta)=\bar{e}_2(\bar{e}_1(\eta))$.
By this way, we show that any point $\eta$ of $A$ admits a finite or an infinite expansion $e_1e_2\cdots $, $e_i\in \overline{\Gamma}$ defined by $f^n(\eta) \in A(e_n)$, $n\geq 1$ where the sequence terminates at $e_n$ if and only if $f^{n}(\eta)\in A$ and $f^{n+1}(\eta)\notin A$.

The sequence is infinite if and only if $\eta\in L_\Gamma$.
Conversely, for any infinite reduced word $e_1e_2\cdots $, the set $\ds\bigcap_{n=1}^{\infty}A(e_1e_2\cdots e_n)=\bigcap_{n=1}^{\infty}f^{-n-1}A(e_n)$ is non-empty and it is reduced to a point belonging on the limit set because the diameter of the sequence of sets $A(e_1e_2\cdots e_n)$ goes to $0$ whenever $n$ goes to $\infty$.

Denote by $\Sigma^{+}$ the set of all infinite reduced words in $\Gamma$. We have showed the following result:

\bpro\label{codlimset}

The map 
\begin{align*}
p^{+}:\Sigma^+&\longrightarrow L_\Gamma\\
                     e_1e_2...&\longmapsto \underset{n\rightarrow+\infty}{\lim} e_1e_2...e_n(o).
\end{align*}
 is a bijection.

\epro
This proposition enables us to represent a point of the limit set $L_\Gamma$ as an infinite reduced word. We will write $\eta=e_1e_2\cdots $ when the infinite reduced word $e_1e_2\cdots $ corresponds to the limit point $\eta$.

\brmq \label{rmk1}

Let $\eta=e_1e_2\cdots $ be a point of the limit set then \\ $\bar{ e}_1(\eta)=e_2e_3\cdots $ and for any $f\in \overline{\Gamma}\setminus \{\bar{ e}_1\}$ we have $f(\eta)=fe_1e_2\cdots.$ 
\ermq

\subsection{Representation of geodesics}
In this section, we use the coding of the limit set in order to give a representation of the oriented closed geodesics of the pair of pants. We begin with the geodesics of the Poincaré disc $\D$. Remember that a geodesic in $\D$ can be specified by its two endpoints. Recall that the non-wandering set of a flow is the set of points which return infinitely often within bounded distance of some given fixed point. 

Denote by $\mathcal{A}$ be the set of oriented geodesics on $\D$ which intersect the fundamental domain $\mathcal{D}$ and both his endpoints belong on $L_\Gamma$
Let $\gamma$ be a geodesic of $\mathcal{A}$ and denote by $\gamma^{+}$ his positive endpoint and by $\gamma^{-}$ his negative endpoint. By the proposition \ref{codlimset}, there exists two infinite reduced words $ e_1e_2\cdots e_n\cdots $ and $f_1f_2\cdots f_n\cdots $   such that $\gamma^{+}=e_1e_2\cdots e_n\cdots  $ and $\gamma^{-}=f_1f_2\cdots f_n\cdots  $
Because of $\gamma\cap\mathcal{D}\neq \emptyset$, we have $e_1\neq f_1$ and the bi-infinite word $\gamma^{+}*\gamma^{-}=\cdots \bar{f}_n\cdots \bar{f}_2\bar{f}_1e_1e_2\cdots e_n\cdots $ is reduced.
We have the following result of C. Series \cite{S2}

\bpro
The map $p:\mathcal{A}\rightarrow \Sigma$ defined by $p(\gamma)=\gamma^{+}*\gamma^{-}$ for all $\gamma \in \mathcal{A}$ is a bijection.

\epro

The following result tells us that the action of the shift $\sigma$ on $\Sigma$ and the action of $\Gamma$ on $\mathcal{A}$ permutes.

\bpro
Let $\gamma$ be a geodesic of $\mathcal{A}$ and $g$ an isometry in $\Gamma$. There exists an integer $N\in \Z$ such that $p(g(\gamma))=\sigma^Np(\gamma)$.
where $\sigma^N=\sigma\sigma\cdots \sigma$

\epro

\brmq
Let $\gamma$ be an oriented geodesic in $\D$ such that both of his endpoints lie in $L_\Gamma$ but $\gamma$ doesn't intersect the fundamental domain $\mathcal{D}$. If $\gamma^{+}=e_1e_2\cdots $ and $\gamma^{-}=f_1f_2\cdots$ then there exists an integer $m$ such that $e_1e_2\cdots e_m=f_1f_2\cdots f_m$ and $e_{m+1}\neq f_{m+1}$. Thus the geodesic $\gamma$ is map by the isometry $g=\bar{e}_m\bar{ e}_{m-1}\cdots \bar{e}_1$ to an oriented geodesic which intersects the fundamental domain $\mathcal{D}$ and both of his endpoints lie in $L_\Gamma$.
\ermq

The two last results tell us that there is a bijection between the set of all orbits of the action of $\sigma$ on the set of bi-infinite reduced words, $\Sigma$ and the set of orbits of the action of $\Gamma$ on the set of all oriented geodesics in $\D$ which intersect the fundamental domain and both of whose endpoints lie in $L_\Gamma$ up to cyclic permutation.
In other words, $\Sigma$ is in bijection with the non-wandering set of the geodesic flow on $P$ up to cyclic permutation .

Now we will focus on the closed geodesics and on the bi-infinite reduced periodic words of $\Sigma$. A word {\bf$W$} is said to be periodic if there exits an integer $n$ such that $\sigma^{n}$({\bf $W$}$)=W$. Thus if {\bf${ W}$} is periodic, there exists a cyclically reduced word $w=e_1e_2\cdots e_n$ such that \\{\bf${W}$}$=\cdots www\cdots $. In this case we will write {\bf${ W}$}$=\langle w\rangle=\langle e_1e_2\cdots e_n\rangle$.

Let $\gamma$ be an oriented geodesic of $P$ in the non-wandering set of the geodesic flow. Assume that the word {\bf${W}$}$=\langle e_1e_2\cdots e_n \rangle$ associated to $\gamma$ is periodic. For any integer $1\leq i\leq n$, denote by $\gamma_i^{+}=e_ie_{i+1}\cdots e_{i-1}$ and by $\gamma_i^{-}=\bar{ e}_{i-1}\bar{ e}_{i-2}\cdots \bar{ e}_i$ the points of $L_\Gamma$ associated to the infinite reduced periodic words $e_ie_{i+1}\cdots $ and $\bar{ e}_{i-1}\bar{e}_{i-2}\cdots$ . Because {\bf${ W}$} is reduced, the geodesic of $\D$ denoted by $\gamma_i$ which endpoints are $\gamma_i^{+}$ and $\gamma_i^{-}$ intersects the fundamental domain $\mathcal{D}$.
\blem\label{lift}
The unique lifts of $\gamma$ in the Poincaré disc $\D$ which intersect the fundamental domain $\mathcal{D}$ are the geodesics $\gamma_i$ for $1\le i \leq n$.

\elem

\bproof
Let $\alpha=(\alpha^{+};\alpha^{-})$ be a lift of $\gamma$ in $\D$ which intersect the fundamental domain $\mathcal{D}$ , then there exists an isometry $g=f_1f_2\cdots f_p \in \Gamma$ where $f_i\in \overline{\Gamma}$ such that $g(\gamma_1)=\alpha$. Thus, by the remark \ref{rmk1}, we have $\alpha^{+}=f_1f_2\cdots f_pe_1e_2\cdots e_ne_1e_2\cdots e_n\cdots $ and $\alpha^{-}=f_1f_2\cdots f_p\bar{e}_n\cdots \bar{ e}_1\cdots .$
Without loss of generality, we can suppose that $p<n$. Then the geodesic $\alpha$ intersects the fundamental domain $\mathcal{D}$ if and only if\\ $f_1f_2\cdots f_p=(e_1e_2\cdots e_p)^{-1}$ or $f_1f_2\cdots f_p=(\bar{e}_n\bar{e}_{n-1}\cdots \bar{e}_{n-p})^{-1}$. In the first case , $\alpha=\gamma_p$ and in the second $\alpha=\gamma_{n-p}$.

\eproof

Together with the fact that any non-wandering geodesic of $P$ is closed if and only if it admits a finite number of lifts in $\D$ which intersect the fundamental domain $\mathcal{D}$, this lemma implies the following result:

\bpro
Any closed geodesic of $P$ is associated to a unique finite cyclically reduced word up to cyclic permutation.

\epro

We will use this symbolic point of view in section \ref{calcul} to calculate the number of self-intersections of the closed geodesics of $P$.

\section{Self-intersection of closed geodesic}\label{calcul}

In this section, we will give a method which enables us to determine the number of self-intersections of a closed geodesic  of $P$ by using the coding  of closed geodesics of $P$ introduced in section \ref{codgeo}.

\subsection{Cyclically lexicographical ordering}

Let $\gamma_1=(\gamma_1^{+};\gamma_1^{-})$  and $\gamma_2=(\gamma_2^{+};\gamma_2^{-})$ be two oriented geodesics in $\D$. The geodesic $\gamma_1$ intersects the geodesic $\gamma_2$ if and only if the endpoints of $\gamma_1$, $\gamma_1^{+}$ and $\gamma_1^{-}$ separates the endpoints of $\gamma_2$, $\gamma_2^{+}$ and $\gamma_2^{-}$.

A {\bf cyclic alphabet} is a cyclically ordered set of distinct symbols and an {\bf alphabet} is a finite ordered set of distinct symbols.\\
$\overline{\Gamma}=\{a, \bar{b}, b, \bar{a}\}$ is the cyclic alphabet whose letters are the generating set arranging in the order in which the first orders intervals (see figure \ref{coding}) occur around $\partial\D$ anticlockwise.
Thus to the cyclic alphabet\\ $\overline{\Gamma}=\{a,\bar{b}, b, \bar{a}\}$ we can associate four distinct alphabets $\Gamma_a=\{a, \bar{a}, b, \bar{a}\}$, $\Gamma_{\bar{b}}=\{\bar{b}, b, \bar{a}, a\}$, $\Gamma_b=\{b, \bar{a}, a, \bar{b}\}$ and $\Gamma_{\bar{a}}=\{\bar{a}, a, \bar{b}, b\}$.

In \cite{BS}, J. Birman and C.Series establish the next theorem:

\bthm[Thm A, \cite{BS}]
Let $\eta=e_1e_2\cdots $ and $\zeta=f_1f_2\cdots $ be two distinct points on $\partial\D$. Then $\eta$ precedes $\zeta$ in anticlockwise order around $\partial\D$ starting from the point $I$ ( see figure \ref{coding}) if and only if:
\benum
\item $e_1$ precedes $f_1$ in the alphabet $\Gamma_a$, or

\item $e_i=f_i$ for $i=1, \cdots , m$ and $e_{m+1}$ precedes $f_{m+1}$ in the alphabet $\Gamma_{\bar{e}_m}$.
\eenum

\ethm
This theorem enables us to define an order on the set of infinite reduced words which is compatible with the representation of points of $L_\Gamma$. From now to see if an oriented geodesic in $\D$ intersect another oriented geodesic, it suffices to look the infinite reduced words associated to the endpoints of these geodesics.

\subsection{Self-intersection and coding}
Let $\gamma$ be a closed geodesic on $P$, it admits a self-intersection on $P$ if there exists two lifts of $\gamma$ on $\D$ which intersect each other in $\D$. If one of those lifts belongs entirely on $D(e)$ for $e\in \overline{\Gamma}$, there exists an isometry $g \in \Gamma$ which maps these two geodesics on two lifts of $\gamma$ on $\D$ which intersect the fundamental domain $\mathcal{D}$ and intersect each other. Thus in order to calculate the number of self-intersections of a closed geodesic, it suffices to consider the lifts which intersect the fundamental domain $\mathcal{d}$. The lemma \ref{lift} enables us to describe all of these lifts if we know the word associated to the closed geodesic $\gamma$.

Let $w=e_1e_2\cdots e_n$ be a  cyclically reduced word. Let $\gamma=\gamma(w)$ be the closed geodesic of $P$ associated to this word. The lifts of $\gamma$ in $\D$ which intersect the fundamental domain $\mathcal{D}$ are the oriented geodesics $\gamma_i$ of positive endpoint $\gamma_i^{+}=e_ie_{i+1}\cdots e_{i-1}$ and negative endpoint $\gamma_i^{-}=\bar{e}_{i-1}\cdots \bar{e}_1\bar{e}_n\cdots \bar{e}_i$ for all $1\leq i\leq n$.
If two geodesics of $\D$, $\gamma_i$ and $\gamma_j$ intersect each other we will note by $(\gamma_i;\gamma_j)$ the self-intersection point.

Denote by $A_\gamma=\{(\gamma_i;\gamma_j); \gamma_i\cap\gamma_j\neq \emptyset; 1\leq i\neq j\leq n\}$ the set of intersection points between two lifts of $\gamma$ which intersect the fundamental domain $\mathcal{D}$. To calculate the number of points of self-intersections of the closed geodesic $\gamma$, it suffices to determine the number of points of $A_\gamma$ lying in the fundamental domain. But it is not always easy to know if two lifts intersect each other in the fundamental or no.

\begin{figure}
 \begin{center}
\includegraphics[scale=0.2]{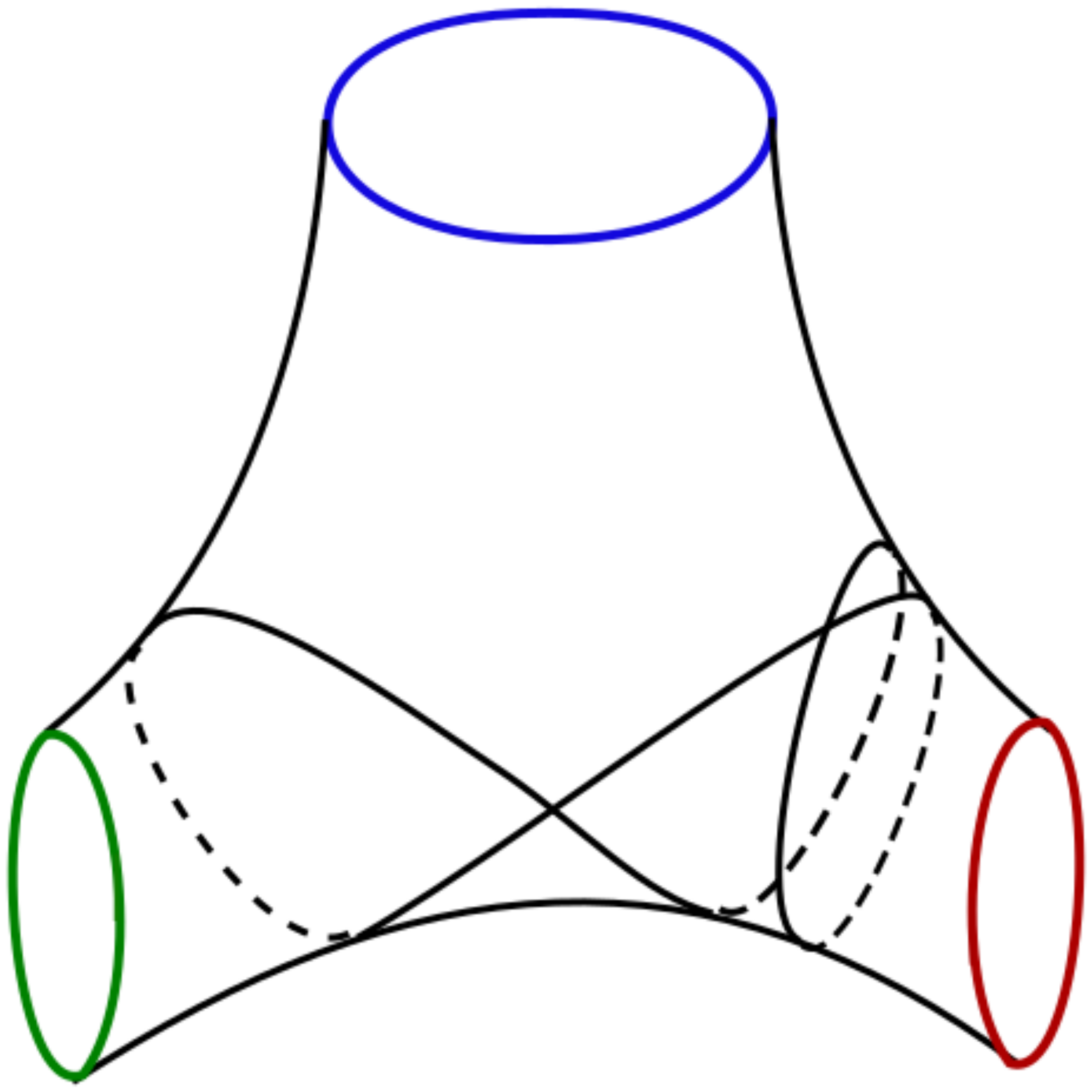}\includegraphics[scale=0.2]{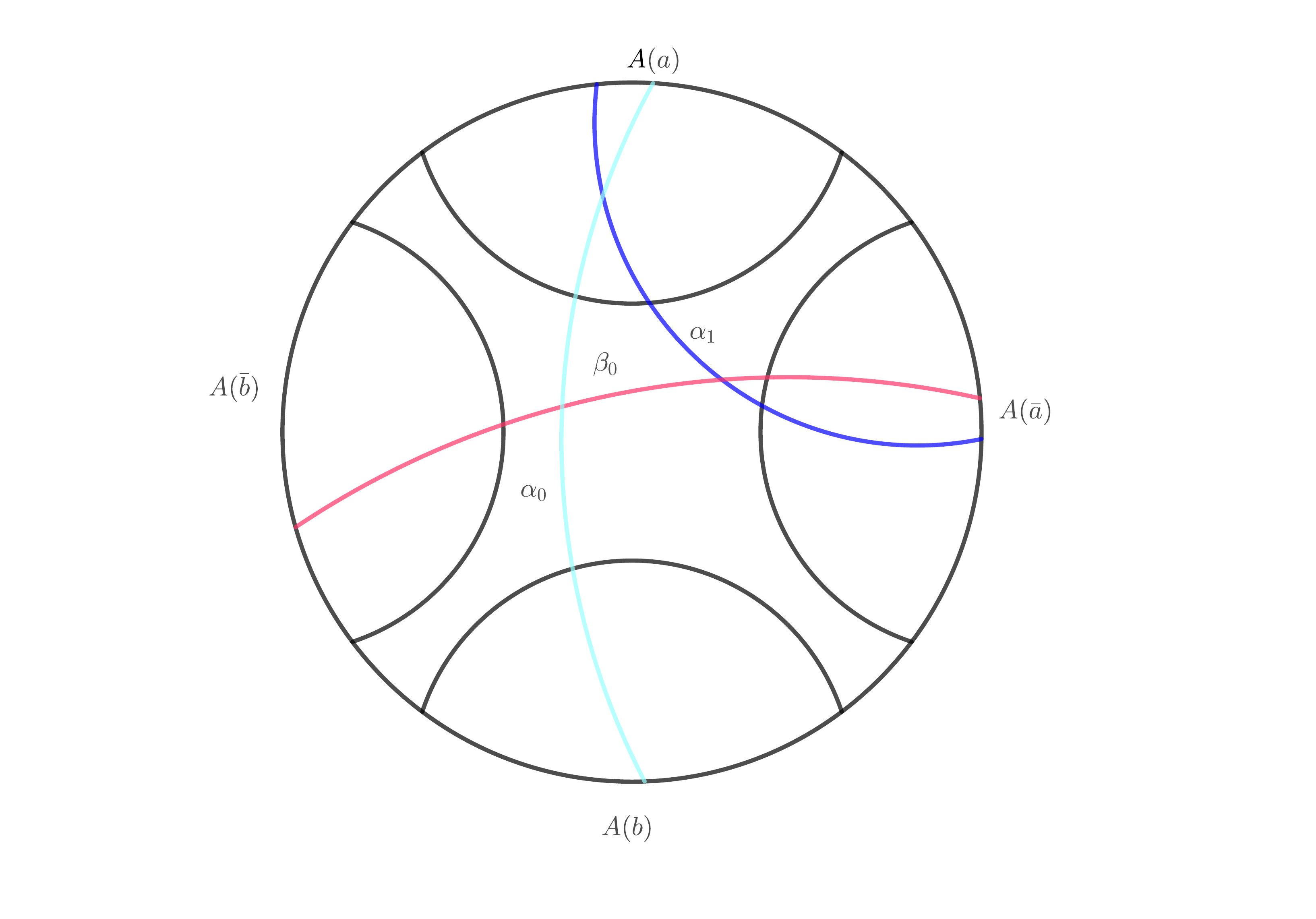}
\caption{}\label{exemple}
  \end{center}
  \end{figure}

For example the closed geodesic $\gamma$ associated to the cyclically reduced non-periodic word $w=a^3\bar{b}$ has two self-intersections (see figure \ref{exemple}). It admits three lifts whose intersect the fundamental domain (see figure \ref{exemple}) $\alpha_0(\alpha_0^{+}=a^2\bar{b};\alpha_0^{-}=b\bar{a}^2)$ , $\alpha_1(\alpha_1^{+}=a\bar{b}a;\alpha_1^{+}=\bar{a}b\bar{a})$ and $\beta_0(\beta_0^{+}=\bar{b}a^2;\beta_0^{-}=\bar{a}^2b)$. The geodesic $\alpha_0$ intersects the geodesic $\beta_0$ and the intersection point lies in the fundamental domain $\mathcal{D}$. Using the cyclically lexicographical ordering, we see that the geodesic $\alpha_1$ intersects the geodesics $\alpha_0$ and $\beta_0$. But we can not say if the points $(\alpha_0;\alpha_1)$ and $(\beta_0;\alpha_1)$ belongs in the fundamental domain. It is easy to see that the isometry $a$ maps the geodesics $\beta_0$ and $\alpha_1$ on the geodesics $\alpha_1$ and $\alpha_0$ and so the the points  $(\alpha_0;\alpha_1)$ and $(\beta_0;\alpha_1)$ are projected on the same point on $P$ and this shows that the number of self-intersections of $\gamma$ is actually $two$.

In order to solve this problem, we define the following equivalence relation on $A_\gamma$: we will say that two points $(\gamma_i;\gamma_j)$ and $(\gamma_k;\gamma_l)$ are equivalent and we write $(\gamma_i;\gamma_j)\sim (\gamma_k;\gamma_l)$ if there exists an isometry in $\Gamma$ which maps the two geodesics $\gamma_i$ and $\gamma_j$ on the geodesics $\gamma_k$ and $\gamma_l$.

{\it The number of self-intersections of closed geodesic $\gamma$ is equal to the cardinal of the quotient of $A_\gamma$ by this equivalence relation:}
$$i(\gamma;\gamma)=\#A_\gamma/\sim.$$

\brmq\label{periodic}
Let $\gamma$ and $\gamma_0$ be two closed geodesics and $m$ an integer such that $\gamma=(\gamma_0)^m$, then: $$i(\gamma;\gamma)=m^2i(\gamma_0;\gamma_0).$$

\ermq

\subsection{Calculation of the number of self-intersections of a closed geodesic}
A cyclically reduced word in $\Gamma$ has the form $w=s_1^{i_1}r_{1}^{j_1}\cdots s_{n}^{j_n}r_n^{j_n}$ where $n, i_1, \cdots i_n,j_1,\cdots, j_n$ are non-zero positive integers and $s_k\in \{a,\bar{a}\}$ and $r_k\in \{b,\bar{b}\}$ for any $1\leq k\leq n$.
From now, every cyclically reduced word will be write in this form. We consider the following cyclic permutations of the word $w$ for all $1\leq k\leq n$, $0\leq m_k\leq i_k-1$ and $0\leq p_k\leq j_k-1$:
$$w_{k,m_k}=s_k^{i_k-m_k}r_k^{j_k}\cdots s_n^{i_n}r_n^{j_n}\cdots s_{k-1}^{i_{k-1}}r_{k-1}^{j_{k-1}}s_k^{m_k}$$
 and
  $$v_{k,p_k}=r_k^{j_k-p_k}s_{k+1}^{j_{k+1}}\cdots r_n^{j_n}s_1^{i_1}\cdots s_k^{i_k}r_k^{p_k}.$$
For $1\leq k\leq n$, we denote respectively by $\alpha_{k,m_k}^{+}$ and $\alpha_{k,m_k}^{-}$ the points of $L_\Gamma$ associated to the infinite reduced periodic words of period  $w_{k,m_k}$ and $w_{k,m_k}^{-1}$ and by $\beta_{k,p_k}^{+}$ and $\beta_{k,p_k}^{-}$ the points of $L_\Gamma$ associated respectively to the infinite reduced periodic words of period $v_{k,p_k}$ et $v_{k,p_k}^{-1}$. The lifts of the closed geodesic $\gamma$ which intersect the fundamental domain  are the geodesics: $\alpha_{k,m_k}=(\alpha_{k,m_k}^{+};\alpha_{k,m_k}^{-})$ and $\beta_{k,p_k}=(\beta_{k,p_k}^{+};\beta_{k,p_k}^{-})$ for all $1\leq k\leq n$, $0\leq m_k\leq i_k-1$ and $0\leq p_k\leq j_k-1$.

We will establish some results which enable us to understand better the intersections between these lifts.
In the next lemma, we show that for every integer $1\leq k\leq n$, the geodesics $\alpha_{k,m_k}$ and $\alpha_{k,m_k^{\prime}}$ for $0\leq m_k<m_k^{\prime}<i_k$ intersect each other  and the intersections points of type $(\alpha_{k,m_k};\alpha_{k,m_k^{\prime}})$ are equivalent to intersection points of type $(\alpha_{k,0};\alpha_{k,m_k^{"}})$ with $0\leq m_k^{"}<i_k$. We show also the same thing for the intersection points of type $(\beta_{k,p_k};\beta_{k,p_k^{\prime}})$ for $0\leq p_k<p_k^{\prime}\leq j_k-1$.
\blem \label{lemint1}
Let $k$ be an integer between  $1$ and $n$.
\benum
\item For any integers $0\leq m_k<m'_k\leq i_k-1,\ \alpha_{k,m_k}\cap \alpha_{k,m'_k}\neq\emptyset$ and we have: 
$$(\alpha_{k,m_k};\alpha_{k,m'_k})\sim(\alpha_{k,0};\alpha_{k,m'_k-m_k})\sim (\beta_{k,0};\alpha_{k,i_k+m_k-m'_k}).$$
\item For any integers $0\leq p_k<p'_k\leq j_k-1,\ \beta_{k,p_k}\cap \beta_{k,p'_k}\neq\emptyset$ and we have: 
$$(\beta_{k,p_k};\beta_{k,p'_k})\sim(\beta_{k,0};\beta_{k,p'_k-p_k})\sim (\alpha_{k+1,0};\beta_{k,j_k+p_k-p'_k}).$$
\eenum
\elem
\bproof 
\benum
\item For any integer $0\leq m_k< m'_k\leq i_k-1,$  with the alphabet $A_a$ , we have $\alpha_{k,m_k}^{+}<\alpha_{k,m'_k}^{+}<\alpha^{-}_{k,m_k}<\alpha^{-}_{k,m'_k}$ if $s_k=a$ and $\alpha^{-}_{k,m'_k}<\alpha^{-}_{k,m_k}<\alpha^{+}_{k,m'_k}<\alpha^{+}_{k,m_k}$ if $s_k=\bar{a}$. So in every case the pair of points $(\alpha^{-}_{k,m_k},\alpha^{+}_{k,m_k})$ separates the pair $(\alpha^{-}_{k,m'_k},\alpha^{+}_{k,m'_k})$ and therefore the geodesics $\alpha_{k,m_k}$ and $\alpha_{k,m'_k}$ intersect each other. Moreover we have: $s_k^{m_k}(\alpha_{k,m_k})=\alpha_{k,0}$, $s_k^{m_k}(\alpha_{k,m'_k})=\alpha_{k,m'_k-m_k},$ $s_k^{m'_k-i_k}(\alpha_{k,m'_k})=\beta_{k,0}$ and $s_k^{m'_k-i_k}(\alpha_{k,m_k})=\alpha_{k,i_k+m_k-m'_k}$. 
It follows that : $$(\alpha_{k,m_k};\alpha_{k,m'_k})\sim (\alpha_{k;0};\alpha_{k,m'_k-m_k})\sim (\beta_{k;0};\alpha_{k,i_k+m_k-m'_k}).$$
\item Likewise for any integers  $0\leq p_k< p'_k\leq j_k-1,$ with the alphabet $A_a$, we obtain $\beta^{-}_{k,p_k}<\beta^{-}_{k,p'_k}<\beta^{+}_{k,p_k}<\beta^{+}_{k,p'_k}$  if $r_k=b$ and $\beta^{+}_{k,p'_k}<\beta^{+}_{k,p_k}<\beta^{-}_{k,p'_k}<\beta^{-}_{k,p_k}$ if $r_k=\bar{b}$. Thus the pair of points $(\beta^{-}_{k,p_k},\beta^{+}_{k,p_k})$ separates the pair $(\beta^{-}_{k,p'_k},\beta^{+}_{k,p'_k})$ and therefore the geodesics $\beta_{k,p_k}$ and $\beta_{k,p'_k}$ intersect each other. Moreover we have: $r_k^{p_k}(\beta_{k,p_k})=\beta_{k,0}$ and $r_k^{p_k}(\beta_{k,p'_k})=\beta_{k,p'_k-p_k}$, and also $r_k^{p'_k-j_k}(\beta_{k,p'_k})=\alpha_{{k+1},0}$ and $r_k^{p'_k-j_k}(\beta_{k,p_k})=\beta_{k,j_k+p_k-p'_k}$. It follows that: $(\beta_{k,p_k};\beta_{k,p'_k})\sim (\beta_{k,0};\beta_{k,p'_k-p_k})\sim (\alpha_{{k+1},0};\beta_{k,j_k+p_k-p'_k})$.  
\eenum
\eproof
In the lemma $\ref{lemint1}$, we fix the integer $k$. Now we want to know what happens between intersection points of type $(\alpha_{k,m_k};\alpha_{l,m_l})$ or $(\beta_{k,p_k};\beta_{l,p_l})$ whenever $k\leq l$. The lemma $\ref{lemint2}$ give us similar results to those of lemma \ref{lemint1}.

{\bf Notation:}   Set $\alpha_{k, i_k}=\beta_k$ for all $1\leq k\leq n$ and $\beta_{k,j_k}=\alpha_{k+1,0}$ for all $1\leq k\leq n$ et $\beta_{n,i_n}=\alpha_{1,i_1}.$
\blem \label{lemint2} Let $k$ and $l$ be two integers between $1$ and $n$. 
\benum
\item For any integers $0\leq m_k\leq i_k-1,\ 0\leq m_l\leq i_l-1,$ \\if $\alpha_{k,m_k}\cap\alpha_{l,m_l}\neq\emptyset,$ then we have :
\pn $\bullet$ if $s_k=s_l,$\\
 \(\left \{  \begin{array}{ll}

(\alpha_{k,m_k};\alpha_{l,m_l})\sim (\alpha_{k;0};\alpha_{l,m_l-m_k})\sim (\beta_{l;0};\alpha_{k,i_l+m_k-m_l}) \\

 \mbox{or}\\
 (\alpha_{k,m_k};\alpha_{l,m_l})\sim (\alpha_{l,0};\alpha_{k,m_k-m_l})\sim (\beta_{k,0};\alpha_{l,i_k+m_l-m_k})
\end{array}\right. \)  
\pn $\bullet$ if $s_k=\overline s_l,$\\
 \(\left \{ \begin{array}{ll}
 (\alpha_{k,m_k};\alpha_{l,m_l})\sim (\alpha_{l,0};\alpha_{k,m_k+m_l})\sim (\alpha_{k;0};\alpha_{l,m_l+m_k})\\

            \mbox{or}\ \\
  (\alpha_{k,m_k};\alpha_{l,m_l})\sim (\beta_{k,0};\alpha_{l,m_l+m_k-i_k})\sim (\beta_{l,0};\alpha_{k,m_k+m_l-i_l})
\end{array} \right. \)
\md
\item For any integers $0\leq p_k\leq j_k-1,\ 0\leq p_l\leq j_l-1,$ \\ if $\beta_{k,p_k}\cap\beta_{l,p_l}\neq\emptyset,$ then we have :
\pn $\bullet$ if $r_k=r_l,$ \\ \(\left \{ \begin{array}{ll}
(\beta_{k,p_k};\beta_{l,p_l})\sim (\beta_{k,0};\beta_{l,p_l-p_k})\sim (\alpha_{l,0};\beta_{k,j_l+p_k-p_l})\\
 \mbox{or}\ \\
  (\beta_{k,p_k};\beta_{l,p_l})\sim (\beta_{l,0};\beta_{k,p_k-p_l})\sim (\alpha_{k,0};\beta_{l,j_k+p_l-p_k})
\end{array} \right. \)

\pn $\bullet$ if $r_k=\overline r_l,$\\
\(\left \{ \begin{array}{ll}
 (\beta_{k,p_k};\beta_{l,p_l})\sim (\beta_{k,0};\beta_{l,p_l+p_k})\sim (\beta_{l,0};\beta_{k,p_k+p_l})\\ \mbox{or}\ \\
  (\beta_{k,p_k};\beta_{l,p_l})\sim (\alpha_{l,0};\beta_{k,p_k+p_l-j_l})\sim (\alpha_{k,0};\beta_{l,p_l+p_k-j_k})
  \end{array} \right. \)
\md

\eenum
\elem
\bproof

\benum
\item We assume that $s_k=s_l$. If the geodesics $\alpha_{k,m_k}$ and $\alpha_{l,m_l}$ intersect each other, then we have: 
\begin{itemize}
\item either $i_k-m_k\leq i_l-m_l$ and $m_l\leq m_k$
\item  $i_k-m_k\geq i_l-m_l$ and $m_l\geq m_k$
\end{itemize}
In the first case the $s_k^{m_k}$ maps the geodesics $\alpha_{k,m_k}$ and $\alpha_{l,m_l}$ respectively on the geodesics $ \alpha_{k,0}$ and $\alpha_{l,m_l-m_k}$ and the isometry $s_k^{m_l-i_l}$ maps the geodesics $\alpha_{k,m_k}$ and $\alpha_{l,m_l}$ respectively on the geodesics  $\alpha_{k,i_l+m_k-m_l}$ et $\beta_{l;0}$. Therefore the intersection point $(\alpha_{k,m_k};\alpha_{l,m_l})$ is equivalent to the intersection points $ (\alpha_{k;0};\alpha_{l,m_l-m_k})$ and $ (\beta_{l,0};\alpha_{k,i_l+m_k-m_l})$. In the second case the isometries  $s_k^{m_l}$ and $\bar{s}_k^{i_k-m_k}$ enable us to have:$$(\alpha_{k,m_k};\alpha_{l,m_l})\sim (\alpha_{l,0};\alpha_{k,m_k-m_l})\sim (\beta_{k,0};\alpha_{l,i_k+m_l-m_k}).$$
Likewise, if $s_k=\bar{s}_l$ and if $\alpha_{k,m_k}\cap\alpha_{l,m_l}\neq \emptyset,$ we have two cases:
\begin{itemize}
\item $i_k-m_k\geq m_l$ and $i_l-m_l\geq m_k$ 
\item $i_k-m_k\leq m_l$ and $i_l-m_l\leq m_k$.
\end{itemize}
In the first case the isometries  $s_k^{m_k}$ and $s_l^{m_l}$ give us the following relations: $(\alpha_{k,m_k};\alpha_{l,m_l})\sim (\alpha_{k,0};\alpha_{l,m_l+m_k})\sim (\alpha_{l,0};\alpha_{k,m_k+m_l})$. 
In the second case we have:  $$(\alpha_{k,m_k};\alpha_{l,m_l})\sim (\beta_{l,0};\alpha_{k,m_k+m_l-i_l})\sim (\beta_{k,0};\alpha_{l,m_l+m_k-i_k})$$ because of the isometries $s_k^{m_k-i_k}$ and $s_l^{m_l-i_l}$.

\item To prove this part, it suffices to replace $s_k,s_l, \alpha_{k,m_k}$ and  $\alpha_{l,m_l}$ respectively by $r_k,r_l, \beta_{k,p_k}$ and  $\beta_{l,p_l}$ and follow the same way like the point $1$.

\eenum
\eproof

\brmq \label{rmqeq}
\benum
Let $k$ and $l$ be two integers in $[1;n]$;
\item Let $m_k$ and $m_l$ such that $\alpha_{k,m_k}\cap\alpha_{l,m_l}\neq\emptyset.$
\begin{itemize}
\item If $s_k=s_l$ and $m_l-m_k=i_l-i_k$ then $$(\alpha_{k,0};\alpha_{l,i_l-i_k})\sim(\beta_{k,0};\beta_{l,0}).$$

Furthermore if $i_k=i_l$ then $(\alpha_{k,0};\alpha_{l,0})\sim(\beta_{k,0};\beta_{l,0})$.

\item If $s_k=\bar{s}_l$ and $m_l+m_k=i_k$ then $$(\alpha_{k,0};\alpha_{l,i_k})\sim(\beta_{k,0};\alpha_{l,0}).$$

Furthermore if $i_k=i_l$ then $(\beta_{k,0};\alpha_{l,0})\sim(\alpha_{k,0};\beta_{l,0})$.
\end{itemize}

\item Let $p_k$ and $p_l$ such that $\beta_{k;p_k}\cap\beta_{l;p_l}\neq \emptyset$.

\begin{itemize}
\item If $r_k=r_l$ and $p_l-p_k=j_l-j_k$ then $$(\beta_{k,0};\beta_{l,j_l-j_k})\sim(\alpha_{k+1,0};\alpha_{l+1;0}).$$

Furthermore if $j_k=j_l$ then $(\beta_{k,0};\beta_{l,0})\sim(\alpha_{k+1,0};\alpha_{l+1;0})$.

\item If $r_k=\bar{r}_l$ and $p_l+p_k=_k$ then $$(\beta_{k,0};\beta_{l,j_k})\sim(\alpha_{k+1,0};\beta_{l,0}).$$

Furthermore if $j_k=j_l$ then  $(\alpha_{k+1,0};\beta_{l,0})\sim(\beta_{k,0};\alpha_{l+1,0})$
\end{itemize}

\item For any integers $1\leq m_k\leq i_k-1$ and $1\leq p_k\leq j_l-1$ we have:
$$\alpha_{k,m_k}\cap\beta_{l,p_l}=\emptyset.$$
\eenum
\ermq

For any integers $k$ and $l$ between $1$ and $n$:
\begin{itemize}

\item[-]  If $s_k=s_l$, then $\alpha_{k,0} \cap \alpha_{l,m_l}\neq \emptyset$ if $i_k>i_l-m_l$.

\item[-]   If $s_k=\bar{s}_l$, then $\alpha_{k,0} \cap \alpha_{l;m_l}\neq \emptyset$ if $i_k>m_l$
and $\beta_{k,0}\cap \alpha_{l,m_l}\neq \emptyset$\\ if $i_k>i_l-m_l$. 
\end{itemize}
This remark , the lemmas \ref{lemint1} and \ref{lemint2} and the remark \ref{rmqeq} motivates the definition of the following sets for any integer $1\leq k\leq n$:

$$A_k^1=\{(\alpha_{k,0};\ \alpha_{l,m_l});\ \  1\leq l\leq n;\ \  s_k=s_l;\ \ \max(0;i_l-i_k)<m_l<i_l\};$$
$$A_k^2=\{(\alpha_{k,0};\ \alpha_{l,m_l});\ \ k< l\leq n;\ \ s_k=\bar{s}_l;\ \ 0<m_l<\max(i_k;i_l)\};$$
$$A_k^3=\{(\beta_{k,0};\ \alpha_{l,m_l});\ \  k< l\leq n;\ \ s_k=\bar{s}_l;\ \ \max(0;i_l-i_k)<m_l<i_l\}.$$

Likewise, considering geodesics $\beta_{k,0}$, $\alpha_{k+1,0}$ and $\beta_{l,p_l}$, we have the same situation and we define the following sets:

$$B_k^1=\{(\beta_{k,0};\ \beta_{l,p_l});\ \ 1\leq l\leq n;\ \ r_k=r_l;\ \ \max(0;j_l-j_k)<p_l<j_l\};$$
$$B_k^2=\{(\beta_{k,0};\ \beta_{l,p_l});\ \ k< l\leq n;\ \ r_k=\bar{r}_l;\ \ 0<p_l<\max(j_k;j_l)\};$$
$$B_k^3=\{(\alpha_{k+1,0};\ \beta_{l,p_l});\ \ k< l\leq n;\ \ r_k=\bar{r}_l;\ \ \max(0;j_l-j_k)<p_l<j_l\}.$$

\bpro\label{propint1}
Any two points of the set $\ds\bigcup_{i=1}^{3}\bigcup_{k=1}^{n}A_k^i$ are not equivalent.
\epro
\bproof
We shall first prove that two points $(\alpha_{k,0};\alpha_{l,m_l})$ and $(\alpha_{t,0};\alpha_{q,m_q})$ of the set $\ds\bigcup_{k=1}^{n}A_k^1\cup A_k^2$ are not equivalent. We assume that there exists an isometry $g\in \Gamma$ which maps the geodesics $\alpha_{k,0}$ and $\alpha_{l,m_l}$ on the geodesics $\alpha_{t,0}$ and $\alpha_{q,m_q}$.
Without loss of generality, we can suppose that $g(\alpha_{k,0})=\alpha_{q,m_q}$ and $g(\alpha_{l,m_l})=\alpha_{t,0},$ because the other can be solved with the same arguments.
Letting

$f=\left\{
   \begin{array}{cc}
     s_{q}^{m_q}  &  \mbox{if}\  k=q \\
    s_{k}^{i_k} r_{k}^{j_k}\cdots r_{q-1}^{j_{q-1}}s_{q}^{m_q}  &\mbox{if}\ k\neq q
    \end{array}
  \right.$, we obtain $f^{-1}(\alpha_{k,0})=\alpha_{q,m_q}$ and $h^{-1}(\alpha_{l,m_l})=\alpha_{t,0}$ for 
  $ \ h=\left\{
   \begin{array}{cc}
     \overline{s}_{l}^{m_l}  & \mbox{if}\ l=t \\
    s_{l}^{i_l-m_l} r_{l}^{j_l}\cdots s_{t-1}^{i_{t-1}}r_{t-1}^{j_{t-1}} &\mbox{if}\ l\neq t
    \end{array}
  \right.$. This implies that $fg(\alpha_{k,0})=\alpha_{k,0}$ and $hg(\alpha_{l,m_l})=\alpha_{l,m_l}$. As the isometries $w_k$ and $w_{l,m_l}$ fix respectively the geodesics $\alpha_{k,0}$ and $\alpha_{l,m_l}$, then there exists two integers $u$ and $v$ such that: $$g=f^{-1}(w_k)^{u}=h^{-1}(w_{l,m_l})^{v}.$$ Thus g is represented by two different reduced words. It is absurd because $\Gamma$ is a free group.
  
It remains us to show that any point of the set $\ds\bigcup_{k=1}^{n}A_k^1\cup A_k^2$ is not equivalent to a point of $\ds\bigcup_{k=1}^{n}A_k^3$.
Let $(\alpha_{k,0};\alpha_{l,m_l})\in  \ds\bigcup_{k=1}^{n}A_k^1\cup A_k^2$ and $(\beta_{t,0};\alpha_{q,m_q})\in \ds\bigcup_{k=1}^{n}A_k^3$ be two points and assume that there exists an isometry $g\in \Gamma$ such that $g(\alpha_{k,0})=\beta_{t,0}$ and $g(\alpha_{l,m_l})=\alpha_{q,m_q}$.
  
We have $f(\alpha_{q,m_q})=\alpha_{k,0}$ and $h(\alpha_{l,m_l})=\beta_{t,0}$ with
$$f=\left\{\begin{array}{ll}
       \ \ \ \ \ \ \ \  s_{q}^{m_q}  &  \mbox{if}\  k=q \\
       s_{k}^{i_k} r_{k}^{j_k}\cdots r_{q-1}^{j_{q-1}}s_{q}^{m_q}  &\mbox{if}\ k\neq q
       \end{array}
     \right.$$ and $$ h=\left\{
      \begin{array}{ll}
       \ \ \ \ \ \ \ \  s_{l}^{i_l-m_l}  & \mbox{if}\ l=t \\
       s_{l}^{i_l-m_l} r_{l}^{j_l}\cdots r_{t-1}^{j_{t-1}}s_{t}^{i_t} &\mbox{if}\ l\neq t
       
      \end{array}
     \right..$$
Then the isometries $fg$ and $hg$ fix respectively the geodesics $\alpha_{k,0}$ and $\alpha_{l,m_l}$. So there exists two integers $u$ and $v$ such that $$g=f^{-1}(w_k)^{u}=h^{-1}(w_{l;m_l})^{v}.$$ As the precedent point, $g$ is represented by two different reduced words in a free group. 

\eproof
By replacing $\alpha_{k,0}$ by $\beta_{k,0}$ and $\alpha_{l,m_l}$ by $\beta_{l,p_l}$ we have:
\bpro\label{propint2}
Any two points of the set $\ds\bigcup_{i=1}^{3}\bigcup_{k=1}^{n}B_k^i$ are not equivalent.
\epro
With the same arguments we prove the following result.
\bpro\label{propint3}
The points of $\ds\bigcup_{i=1}^{3}\bigcup_{k=1}^{n}A_k^i$ are not equivalent to the points of $\ds\bigcup_{i=1}^{3}\bigcup_{k=1}^{n}B_k^i$.
\epro

Now, to compute the number of self-intersections, it remains us to examine the intersection points of type $(\alpha_{k,0};\alpha_{l,0})$, $(\alpha_{k,0};\beta_{l,0})$ and $(\beta_{k,0};\beta_{l,0})$.

We will show first that these points are not equivalent to the points of $A_i$ and $B_i$ for $i=1,2,3$.
\bpro
The intersection points of type $(\alpha_{k,0};\alpha_{l,0})$, $(\alpha_{k,0};\beta_{l,0})$ and $(\beta_{k,0};\beta_{l,0})$ are not equivalent to the point of the set $\ds\bigcup_{i=1}^{3}\bigcup_{k=1}^{n}A_k^i\cup B_k^{i}$.
\epro
\bproof
It suffices to take $m_l=0$ or $p_l=0$ in the proofs of propositions \ref{propint1}, \ref{propint2} and \ref{propint3}.
\eproof

The remark \ref{rmqeq} motivates the definition of the following sets for any integer $1\leq k\leq n$:

$$C_k^1(w)=\{(\alpha_{k,0};\alpha_{l,0}); \ \  \alpha_{k,0}\cap\alpha_{l,0}\neq\emptyset;\ \  s_k^{i_k}\neq s_l^{i_l}; \ \ k<l\leq n \};$$
$$C_k^2(w)=\{(\beta_{k, 0};\alpha_{l,0});\ \ \beta_{k,0}\cap\alpha_{l,0}\neq \emptyset; \ s_k^{i_k}\neq \bar{s}_l^{i_l};\ \  k<l\leq n\};$$
$$D_k^1(w)=\{(\beta_{k,0};\beta_{l,0}); \ \ \beta_{k,0}\cap\beta_{l,0}\neq\emptyset;\ \ r_k^{j_k}\neq r_l^{j_l}; \ \ k<l\leq n \};$$
$$D_1^2(w)=\{(\alpha_{1,0};\beta_{l,0}); \ \ \alpha_{1,0}\cap\beta_{l,0}\neq \emptyset; \ \ \ 1\leq l\leq n\};$$
and for all $2\leq k\leq n$,
$$D_k^2(w)=\{(\alpha_{k,0};\beta_{l,0});\ \  \alpha_{k,0}\cap\beta_{l,0}\neq\emptyset;\ \ \ r_{k-1}^{j_{k-1}}\neq \bar{r}_{l}^{j_l}; \ \ k\leq l \leq n\}.$$

\bpro
Let $k$ and $l$ be two integers between $1$ and $n$ such that $k\leq l$ and let $w=s_1^{i_1}r_1^{j_1}\cdots s_n^{i_n}r_n^{j_n}$ a finite non-periodic word, then we have:
\benum
\item 
Every intersection point of type $(\alpha_{k,0};\alpha_{l,0})$ or $(\beta_{k,0};\beta_{l,0})$ is equivalent to a point of $\ds\bigcup_{k=1}^{n}C_k^1(w)\cup D_k^1(w).$ 
\item
Every intersection point of type $(\alpha_{k,0};\beta_{l,0})$ is equivalent to a point of $\ds\bigcup_{k=1}^{n}C_k^2(w)\cup D_k^2(w).$

\item 
Two different points of $\ds\bigcup_{i=1}^{2}\bigcup_{k=1}^{n}C_k^i(w)\cup D_k^i(w)$ are not equivalent.
\eenum
\epro

\bproof
We assume that $\alpha_{k,0}\cap  \alpha_{l,0}\neq \emptyset$. If $s_k^{i_k}\neq s_l^{i_l}$, then $(\alpha_{k,0};\alpha_{l,0})\in\ds \bigcup_{k=1}^{n}C_k^i(w) $ .  In the case where  $s_k^{i_k}=s_l^{i_l}$, then  $\bar{s}_k(\alpha_{k,0})=\beta_{k,0}$ and $\bar{s}_l(\alpha_{l,0})=\beta_{l,0}$. Thus $\beta_{k,0}\cap \beta_{l,0}\neq \emptyset$ and the points $(\alpha_{k,0};\alpha_{l,0})$ and $(\beta_{k,0};\beta_{l,0})$ are equivalent. If $r_k^{j_k}\neq r_l^{j_l}$, then $(\beta_{k,0};\beta_{l,0})\in \ds \bigcup_{k=1}^{n}D_k^i(w)$. Otherwise, because $w$ is non-periodic, there exists a positive non-zero integer $p$ such that for any integer $0\leq m< p$ we have:
$$\left\{\begin{array}{lll}
s_{k+m}^{i_{k+m}}=s_{l+m}^{i_{l+m}}\\
r_{k+m}^{j_{k+m}}=r_{l+m}^{j_{l+m}}\\
s_{k+p}^{i_{k+p}}\neq s_{l+p}^{i_{l+p}}\\
\end{array}\right.\ \ \ 
\mbox{or} \ \ \
\left\{\begin{array}{ll}
s_{k+m}^{i_{k+m}}=s_{l+m}^{i_{l+m}}\\
r_{k+m}^{j_{k+m}}=r_{l+m}^{j_{l+m}}\\
s_{k+p}^{i_{k+p}}=s_{l+p}^{i_{l+p}}\\
r_{k+p}^{j_{k+p}}\neq r_{l+p}^{j_{l+p}}\\
\end{array}\right..$$
In the first case we have $g(\alpha_{{k},0})=\alpha_{{k+p},0}$ and 
$g(\alpha_{l,0})=\alpha_{{l+p},0}$ with \\  $g=\bar{r}_{k+p-1}^{j_{k+p-1}}\cdots\bar{r}_{k+m-1}^{j_{k+m-1}}\bar{s}_{k+m-1}^{i_{k+m-1}}\cdots \bar{s}_{i_k}.$
Thus  the points $(\alpha_{k,0}; \alpha_{l,0})$ and $ (\alpha_{{k+p},0}; \alpha_{{l+p},0})$ are equivalent and $ (\alpha_{{k+p},0}; \alpha_{{l+p},0}) \in \ds\bigcup_{k=1}^{n}C_k^1(w) $ . 

 In the second case the point $(\alpha_{k,0};\alpha_{l,0})$ is equivalent to the point $(\beta_{k+p,0};\beta_{l+p,0})$ which lies in $\ds\bigcup_{k=1}^{n}D_k^1(w)$.

We prove the point $2$ with the same arguments.

We prove the point $3$ with the same proof of the proposition \ref{propint1}.
\eproof

Now, we can determine the number of self-intersections $i(\gamma;\gamma)$ of a closed geodesic $\gamma$ of $P$ associated to the cyclically reduced  word  $w=s_1^{i_1}r_1^{j_1}\cdots $ ; set:
$$ H(w)=\ds\sum_{i=1}^{2}\sum_{k=1}^{n}\#C_k^i(w)+\#D_k^i(w).$$
Then the the number of self-intersections of $\gamma$:

$$ i(\gamma;\gamma)=H(w)+\sum_{i=1}^{3}\sum_{k=1}^{n}\#A_k^i+\#B_k^i.$$

\brmq
Let $w$ and $v$ two finite cyclically reduced words and $m$ a non-zero integer such that $w=v^m$. Then we have:

$$ H(w)=m^2H(v).$$

\ermq

Let us compute $\ds\sum_{i=1}^{3}\sum_{k=1}^{n}\#A_k^i+\#B_k^i$ as a function of the integers $n$,  $i_k$ and $j_k$ for $1\leq k\leq n$.
 
 \begin{eqnarray*}
\sum_{k=1}^{n}\#A_k^1 &=& \ds\sum_{k=1}^{n}\sum_{l=1}^{n}\#\{(\alpha_{k,0};\alpha_{l;m_l});\ s_k=s_l;\ \max(0;i_l-i_k)<m_l<i_l \}\\
      &=& \ds\sum_{k=1}^{n}\sum_{l=1}^{n}\big(\min(i_k;i_l)-1\big)\delta_{s_k;s_l}\\
      &=& \sum_{k=1}^{n}i_k-1 \ +\ds 2\sum_{k=1}^{n}\sum_{l=k+1}^{n} \big(\min(i_k;i_l)-1\big)\delta_{s_k;s_l}
  \end{eqnarray*}  
 \begin{eqnarray*} 
 \sum_{k=1}^{n}\#A_k^2 & = & \ds\sum_{k=1}^{n}\sum_{l=1}^{n}\#\{(\alpha_{k,0}; \alpha_{l;m_l}); \  s_k=\bar{s}_l; \ 0<m_l<\min(i_k;i_l) \}\\ 
 & = & \ds \sum_{k=1}^{n}\sum_{l=k+1}^{n}\big(\min(i_k;i_l)-1\big)\delta_{s_k;\bar{s}_l}.
 \end{eqnarray*}
 
 \begin{eqnarray*}
\sum_{k=1}^{n}\#A_k^3 & = & \ds\sum_{k=1}^{n}\sum_{l=1}^{n}\#\{(\beta_{k,0}; \alpha_{l;m_l});\  s_k=\bar{s}_l; \max(0;i_l-i_k)<m_l<i_l\}\\ 
& = &  \sum_{k=1}^{n}\sum_{l=k+1}^{n}\big(\min(i_k; i_l)-1\big)\delta_{s_k;\bar{s}_l}.
\end{eqnarray*} 
 Because  $\ds\delta_{s_k;s_l}+\delta_{s_k;\bar{s}_l}=1$, then we have:
 $$\sum_{i=1}^{3}\sum_{k=1}^{n}\#A_k^i=\sum_{k=1}^{n}i_k-1+ \ 2\sum_{k=1}^{n}\sum_{l=k+1}^{n}\min(i_k;i_l)-1.$$

Similarly, we show that:

 $$\sum_{k=1}^{3}\#B_k=\ds\sum_{k=1}^{n}j_k-1+ \ 2\sum_{k=1}^{n}\sum_{l=k+1}^{n}\min(j_k;j_l)-1.$$

Recall that $L(\gamma)=\ds\sum_{k=1}^{n}i_k+j_k$  is the combinatorial length of $\gamma$.                         
We have just proved the following result:
\bthm\label{autointer}
Let $\gamma$ be a closed geodesic of  $P$ associated to the cyclically reduced word  $w=s_1^{i_1}r_1^{j_1}\cdots s_n^{i_n}r_n^{j_n}$. Then we have:
$$i(\gamma;\gamma)-H(w)=nL(\gamma)-2n^{2}-\ds\sum_{k=1}^{n}\sum_{l=k+1}^{n}|i_k-i_l|+|j_k-j_l|.$$
\ethm
\section{Proofs of theorems \ref{thm1} and \ref{thm2}}

To understand better the self-intersections $i(\gamma;\gamma)$ of a closed geodesic $\gamma$ of $P$, it is necessary to understand $H(w)$.
In this section, we use the theorem \ref{autointer}  in order to prove the theorems \ref{thm1} and \ref{thm2}. 
In this section set: $\alpha_k=\alpha_{k,0}$ and $\beta_k=\beta_{k,0}$. We begin with the following result on $H(w)$ where $w$ is a cyclically reduced non-periodic word.

\bthm\label{H}
Let $\gamma$ be a closed geodesic of $P$ associated to the cyclically reduced word $w(\gamma)=s_1^{i_1}r_1^{j_1}\cdots s_n^{i_n}r_n^{j_n}=(w_0)^q$ where $q\in \N$ and $w_0$ is a  cyclically reduced non-periodic word . Then we have:
$$q(n-q)\leq H(w)\leq n^2+(n-q)^2.$$
In particular if $w$ is cyclically reduced non-periodic word, $$(n-1)\leq H(w)\leq n^2+(n-1)^2.$$
\ethm
We prove first this theorem for the cyclically reduced non-periodic word containing two distinct letters of $\overline{\Gamma}$.
\bpro\label{H1}
Let $\gamma$ be a closed geodesic of $P$ associated to the cyclically reduced non-periodic word $w(\gamma)=s_1^{i_1}r_1^{j_1}\cdots s_n^{i_n}r_n^{j_n}$, then we have:
\benum
\item  If $s_i=a$ and $r_i=\bar{b}$ for all $1\leq i\leq n$, then: $$n^2+n-1\leq H(w)\leq n^2+(n-1)^2.$$

\item If $s_i=a$ and $r_i=b$ for all $1\leq i\leq n$, then: $$n-1\leq H(w)\leq (n-1)^2.$$
\eenum
\epro
\bproof
We consider the cyclically reduced non-periodic word \\$w=a^{i_1}\bar{b}^{j_1}\cdots a^{i_n}\bar{b}^{j_n}$. For all $1\leq k;l\leq n $, the geodesics $\alpha_k$ and $\beta_l$ intersect each other, $s_k=s_l$ and $r_{k-1}=r_l$, then we have:  $$\ds\sum_{k=1}^{n}\#C_k^2(w)+\#D_k^2(w)=n^2.$$ 
It remains us now to prove that: $$n-1\leq \ds\sum_{k=1}^{n}\#C_k^{1}(w)+\#D_k^{1}(w)\leq (n-1)^2.$$
Because $s_k=s_l$ and $r_k=r_l$
for $1\leq k\leq n$,  $$C_k^1(w)=\#\{k<l\leq n:\ \alpha_{k}\cap \alpha_{l}\neq \emptyset \ ; \  \ i_k\neq i_l\}$$
and $$D_k^1(w)=\#\{k<l\leq n: \ \beta_{k}\cap\beta_{l}\neq \emptyset;  \ \ j_k\neq j_l\}.$$

 It follows that:  $$ \ds\sum_{k=1}^{n}\#C_k^1(w)\leq \ds\sum_{k=1}^{n}\#\left\{\begin{array}{ll}
       & i_k< i_l \\
  l:\\
       & j_{k-1}\leq j_{l-1}\\           
      \end{array}\right\}$$
   and   
 $$ \ds\sum_{k=1}^{n}\#D_k^1(w)\leq \ds\sum_{k=1}^{n}
        \#\left\{\begin{array}{ll}
              & i_k\leq i_l \\
         l:\\
              & j_k< j_l\\           
             \end{array}\right\}.$$

There exists a permutation $\sigma$ of $\{i_1,i_2,\cdots , i_n\}$ and a permutation $\tau$ of $\{j_1,j_2,\cdots, j_n\}$ such that $\sigma(i_k)=x_k$ and $\tau(j_k)=y_k$ $\forall$ $1\leq k\leq n$ and for all $1\leq i\leq n-1$, $x_i\leq x_{i+1}$ and $y_i\leq y_{i+1}$.
Thus for any integer $1\leq k\leq n$, there exists three integers $1\leq p\leq n$, $1\leq q\neq m\leq n$ such that $i_k=x_p$, $j_k=y_q$ and $j_{k-1}=y_m$. 
So $$ \ds\sum_{k=1}^{n}\#C_k^1(w)\leq \ds\sum_{k=1}^{n}\#\left\{\begin{array}{ll}
       & x_p<i_l \\
   l:\\
       & y_m\leq j_{l-1}\\
       \end{array}\right\}$$
   and    
 $$\ds\sum_{k=1}^{n}\#D_k^1(w)\leq\ds\sum_{k=1}^{n}\#\left\{\begin{array}{ll}
       & x_p\leq i_l\\
   l: \\
       & y_q<j_l \\
       \end{array}\right\}.$$

We can assume that $m\leq q-1$ and thus: $$ \ds\sum_{k=1}^{n}\#\left\{\begin{array}{ll}
       & x_p<i_l \\
   l:\\
       & y_m\leq j_{l-1}\\
       \end{array}\right\}+\#\left\{\begin{array}{ll}
       & x_p\leq i_l\\
   l: \\
       & y_q<j_l \\
       \end{array}\right\}\leq 2n-2m-1.$$ This implies that: $$\ds\sum_{k=1}^{n}\#C_k^1(w)+\#D_k^1(w)\leq \ds\sum_{m=1}^{n}2n-2m-1=(n-1)^2.$$
      
We use the method to show that the cardinal of the complementary of the set $\ds\bigcup_{k=1}^{n}C_k^1(w)\cup D_k^1(w)$ is less than or equal $(n-1)^2$ and this proves that $\ds\sum_{k=1}^{n}\#C_k^1(w)+\#D_k^1(w)\geq n-1$.

With the same method, we prove that $n-1\leq H(w)\leq (n-1)^2$ for any cyclically reduced non-periodic word of type
 $w=a^{i_1}b^{j_1}\cdots a^{i_1}b^{j_n}$.

\eproof

\brmq

\bitem

\item If $w=a^{i_1}\bar{b}^{j_1}\cdots a^{i_n}\bar{b}^{j_n}$ is periodic then there exists two integers $p, q\in\N^{*}$ such that $w^{\prime}=a^{i_1}\bar{b}^{j_1}\cdots a^{i_p}\bar{b}^{j_p}$ is a non-periodic sub-word of $w$ and $w=(w^{\prime})^q$. By using the proposition \ref{H1} and the fact that $H(w)=q^2H(w^{\prime})$ we have the following inequalities for these words:

$$ n^2+q(n-q)\leq H(w)\leq n^2+(n-q)^2.$$

 \item If we replace $\bar{b}$ by $b$ in the word $w$, we have: $$ q(n-q)\leq H(w)\leq (n-q)^2.$$
\eitem

 \ermq

Now we will prove the theorem \ref{H} for the words with $3$ or $4$ distinct letters of $\overline{\Gamma}$. 
We can suppose that every cyclically reduced word can be written in the following form: $w=a^{i_1}b^{j_1}\bar{a}^{i_2}b^{j_2}\cdots \bar{a}^{i_p}b^{j_p}a^{i_{p+1}}x$ or $w=a^{i_1}b^{j_1}\bar{a}^{i_2}b^{j_2}\cdots \bar{a}^{i_p}b^{j_p}\bar{a}^{i_{p+1}}x$ where $x$ is a sub-word of $w$. We will just deal with the case where $w=a^{i_1}b^{j_1}\bar{a}^{i_2}b^{j_2}\cdots \bar{a}^{i_p}b^{j_p}a^{i_{p+1}}x$ because the other case can be deducted from the first.
Set $u=b^{j_1}\bar{a}^2b^{j_2}\cdots \bar{a}^{i_p}b^{j_p}$ and consider the cyclically reduced word $w^{\prime}=a^{i_1}\bar{u}a^{i_{p+1}}x$. We have:
\bpro
Let $w=a^{i_1}ua^{i_{p+1}}x$ and $w^{\prime}=a^{i_1}\bar{u}a^{i_{p+1}}x$ two cyclically reduced  words, then we have $H(w)\leq H(w^{\prime}).$
\epro
\bproof
Consider the geodesics $\beta_1(\beta_1^{+}=b^{j_1}\cdots xa^{i_1}; \beta_1^{-}=\bar{a}^{i_1}\bar{x}\cdots \bar{b}^{j_1})$ and $\alpha_{p+1}(\alpha_{p+1}^{+}=a^{i_{p+1}}x\cdots b^{j_p};\alpha_{p+1}^{-}=\bar{b}^{j_p}\cdots \bar{x}\bar{a}^{i_{p+1}})$; the geodesic $\beta_1$ doesn't intersect $\alpha_{p+1}$. The transformation introduced changes these two geodesics into two geodesics $\beta_{1}^{\prime}(\beta_{1}^{{\prime}+}=\bar{b}^{j_p}a^{i_p}\cdots xa^{i_1};\beta_{1}^{{\prime}-}=\bar{a}^{i_1}\bar{x}\cdots \bar{a}^{i_p}b^{j_p})$ and $\alpha_{p+1}^{\prime}(\alpha_{p+1}^{\prime}=a^{i_{p+1}}xa^{i_1}\cdots \bar{b}^{j_1};\alpha_{p+1}^{\prime}=b^{j_1}\cdots \bar{x}\bar{a}^{i_{p+1}})$ and these two geodesics intersect each other. By using the definition of $H(w)$, we show that for any integer $1\leq k\leq p$, 
\bitem
\item if $(\alpha_{p-k+2}^{\prime};\alpha_{p+1}^{\prime})\notin C_{p-k+2}^1(w^{\prime})$ then: $$(\beta_k;\alpha_{p+1})\notin C_k^2(w) \ \  \mbox{or}\ \ (\beta_k;\beta_1)\notin D_1^1(w)$$ and 
\item if $(\beta_{p-k+2}^{\prime};\beta_{1}^{\prime})\notin D_1^{1}(w^{\prime})$ then: $$(\beta_p;\beta_{k-1})\notin D_{k-1}^1(w)\ \mbox{or} \ (\alpha_k;\alpha_{p+1})\notin C_{k}^1(w) \ \mbox{or} \ (\beta_1;\alpha_k)\notin C_1^2(w).$$
\eitem

By proceeding in the same way, we also show without difficulty that for $p+1\leq l\leq n$, if $s_{l+1}=a$ and $r_l=\bar{b}$, then if $(\alpha_{l+1}^{\prime};\alpha_{p+1}^{\prime})\notin C_{p+1}^1(w^{\prime})$ then $(\beta_1;\alpha_{l+1})\notin C_1^2(w)$ or $(\alpha_{l+1};\alpha_{p+1})\notin C_{p+1}^1(w)$.

The remaining cases are treated in the same way.
\eproof
The idea is to "straighten" the geodesics of type $(W^{+}=axb;W^{-}=\bar{b} \bar{x}\bar{a})$ and $(V^{+}=bya; V^{-}=\bar{a}\bar{y}\bar{b})$ and to turn them into geodesics of type $(W'^{+}=ax\bar{b}; W'^{-}=b\bar{x}\bar{a})$ or $(V'^{+}=by\bar{a}; V'^{-}=a\bar{y}\bar{b})$.
By repeating this procedure, we get at the end a word containing only the letters $a$ and $\bar{b}$, we obtain:

\bpro\label{H2}
Let $w=s_1^{i_1}r_1^{j_1}\cdots s_n^{i_n}r_n^{j_n}$ a cyclically reduced  word then there exists a permutation $\sigma$ of $\{i_1,i_2,\cdots, i_n\}$ and a permutation $\tau$ of $\{j_1,j_2,\cdots, j_n\}$ and  a cyclically reduced word  \\ $w^{\prime}=a^{\sigma(i_1)}\bar{b}^{\tau(j_1)}\cdots a^{\sigma(i_n)}\bar{b}^{\tau(j_n)}$ such that $H(w)\leq H(w^{\prime})$.

\epro
To prove the other inequality of Theorem \ref{H} for words containing at least 3 letters of $\overline{\Gamma}$, we introduce a transformation that consists of doing the "opposite" of what we did to prove the last proposition \ref{H2}. 
Without loss of generality, we can assume that every word is written as follows: $w=a^{i_1}\bar{b}^{j_1}\bar{a}^{i_2}\cdots \bar{a}^{i_p}\bar{b}^{j_p}a^{i_{p+1}}x$ or $w=a^{i_1}\bar{b}^{j_1}\bar{a}^{i_2}\cdots \bar{a}^{i_p}\bar{b}^{j_p}\bar{a}^{i_{p+1}}x$ where $x$ is a sub-word of $w$. 

Let $v=\bar{b}^{j_1}\bar{a}^{i_2}\cdots \bar{a}^{i_p}\bar{b}^{j_p}$ and $w^{\prime\prime}=a^{i_1}\bar{v}a^{i_{p+1}}x$, the method used in the proof of the proposition \ref{H2} can be used to prove that $H(w^{\prime\prime})\leq H(w)$. 
By repeating this procedure, we get at the end a word containing only the letters $a$ and $b$, we obtain:

\bpro\label{H3}
Let $w=s_1^{i_1}r_1^{j_1}\cdots s_n^{i_n}r_n^{j_n}$ a cyclically reduced  word then there exists a permutation $\sigma$ of $\{i_1,i_2,\cdots, i_n\}$ and a permutation $\tau$ of $\{j_1,j_2,\cdots, j_n\}$ and  a cyclically reduced word\\  $w^{\prime\prime}=a^{\sigma(i_1)}b^{\tau(j_1)}\cdots a^{\sigma(i_n)}b^{\tau(j_n)}$ such that $H(w^{\prime\prime})\leq H(w)$.

\epro

{\bf Proof of theorem  \ref{thm1} }

From theorem \ref{H} and theorem \ref{autointer}, we will deduce the theorem \ref{thm1}.

We begin with the first inequality.

\bpro
Let $\gamma$ be a non-simple closed geodesic of $P$ associated to the cyclically reduced  word $w(\gamma)=s_1^{i_1}r_1^{j_1}\cdots s_n^{i_n}r_n^{j_n}$, then we have: $$i(\gamma;\gamma)\geq L(\gamma)-n-1.$$

\epro

\bproof

Thanks to theorem \ref{autointer}, we have: $$i(\gamma;\gamma)=L-2n+H(w)+2\ds\sum_{1\leq k<l\leq n}\min(i_k;i_l)+\min(j_k;j_l)-2.$$
 If $\ds\sum_{1\leq k<l\leq n}\min(i_k;i_l)+\min(j_k;j_l)-2\neq 0$ then : $$\ds\sum_{1\leq k<l\leq n}\min(i_k;i_l)+\min(j_k;j_l)-2\geq n-1 \  \mbox{and we have the result}.$$

If $\ds\sum_{1\leq k<l\leq n}\min(i_k;i_l)+\min(j_k;j_l)-2=0$, then: $$\#\{1\leq k\leq n: i_k=j_k=n-1\}\geq n-1$$ and we can have two different cases:

\benum

\item  there exists $1\leq k_0\leq n$ such that $i_{k_0}\neq 1$ or $j_{k_0}\neq 1$. In this case, proposition \ref{H} gives us $H(w)\geq n-1$.

\item for any $1\leq k\leq n$, $i_k=j_k=1$, because $\gamma$ is non-simple $w\neq (ab)^n$. Thanks to proposition \ref{H3}, $H(w)\geq \min( H(w_0); H(w_1))$ where $w_1=a(\bar{b}\bar{a})^p\bar{b}(ab)^{n-p-1}$ and $w_1=a(\bar{b}\bar{a})^pb(ab)^{n-p-1}$. It is not difficult to show that $H(w_1)=n-1$ and $H(w_2)=n$.  Thus we have the result.

\eenum

\eproof

Finish the proof of theorem \ref{thm1}.
\bpro
Let $\gamma$ be a closed geodesic of $P$ associated to the cyclically reduced word $w(\gamma)=s_1^{i_1}r_1^{j_1}\cdots s_n^{i_n}r_n^{j_n}$, then we have:

$$ i(\gamma;\gamma)\leq nL(\gamma)-n^2.$$

\epro

\bproof 

As $$i(\gamma;\gamma)=nL(\gamma)-2n^2+H(w)-\ds\sum_{1\leq k<l\leq n}|i_k-i_l|+|j_k-j_l|,$$
it suffices to prove the following inequality: $$H(w)-\ds\sum_{1\leq k<l\leq n}|i_k-i_l|+|j_k-j_l| \leq n^2.$$
By the definition of $H(w)$, we have:
 
   $$H(w)\leq 2n^2-n-\#\{1\leq k<l\leq n; i_k\neq i_l\}-\#\{1\leq k<l\leq n; j_k\neq j_l\}$$
   
   and thus 
   
   $$H(w)\leq 2n^2-n- \ds\sum_{1\leq k<l\leq n}
         \delta_{i_{k};i_{l}}+\delta_{j_{k};j_{l}}.$$
        
  As $\delta_{i,j}= 1$ ou $0$ for any integers $i$ and $j$, we have also: $$\ds\sum_{1\leq k<l\leq n}|i_{k}-i_{l}|+
               \delta_{i_{k};i_{l}}+|j_{k}-j_{l}|+\delta_{j_{k};j_{l}}\geq n^2-n.$$
This implies the result.

\eproof

This completes the proof of the theorem \ref{thm1}.

{\bf Proof of theorem \ref{thm2}}

The theorem \ref{thm2} is deduced from the inequalities of  theorem \ref{thm1}.
\bcor
Let $\gamma$ be a closed geodesic of $P$ and $L(\gamma)$ his  combinatorial length then we have:

$$ i(\gamma;\gamma)\geq \ds\left\{\begin{array}{ll}
\ds\frac{L(\gamma)}{2}-1 \ \  \mbox{if}  \  \mbox{$L(\gamma)$ is even}\\
\\
\ds\frac{L(\gamma)-1}{2} \ \ \mbox{if} \mbox{ $L(\gamma)$ is odd} \\
\end{array}\right.$$

and these bounds are sharp.

The geodesics associated to the following words $w=a(\bar{b}\bar{a})^p\bar{b}(ab)^{n-p-1}$ when $L(\gamma)$ is even and $w=a^{2}b(ab)^{n-1}$ when $L(\gamma)$ is odd realize the minimal self-intersections number.

\ecor

\bproof Consider the closed geodesics $\gamma_1$ and $\gamma_3$ of $P$ associated respectively to the following  cyclically reduced words $w_1=a(\bar{b}\bar{a})^p\bar{b}(ab)^{n-p-1}$ and $w_3=a^{2}b(ab)^{n-1}$ . 
By definition of $H(w)$, we have $$H(w_1)=H(w_3)=n-1.$$ Thus,  by the theorem \ref{autointer}, $$i(\gamma_1;\gamma_1)=i(\gamma_3;\gamma_3)=L(\gamma)-n-1.$$ This gives the proof.

\eproof

\bcor
Let $\gamma$ be a closed geodesic of $P$ and $L(\gamma)$ his  combinatorial length then we have:

$$ i(\gamma;\gamma)\leq  \ds\left\{\begin{array}{ll}
\ds\frac{L^2(\gamma)}{4} \ \ \mbox{if}  \  \mbox{$L(\gamma)$ is even}\\
\\
\ds\frac{L^2(\gamma)-1}{4} \ \ \mbox{if}  \ \mbox{ $L(\gamma)$ is odd} \\
\end{array}\right.$$

and these bounds are sharp.
The geodesics associated to the following words $w=(a\bar{b})^n$ when $L(\gamma)$ is even and $w=a^{2}\bar{b}(a\bar{b})^{n-1}$ when $L(\gamma)$ is odd realize the maximal self-intersections number.

\ecor

\bproof Consider the closed geodesics $\gamma_1^{\prime}$ and $\gamma_3^{\prime}$ of $P$ associated respectively to the following  cyclically reduced words $w_1^{\prime}=(a\bar{b})^n$ and $w_3^{\prime}=a^{2}\bar{b}(a\bar{b})^{n-1}$ . 
By definition of $H(w)$, we have $$H(w_1^{\prime})=n^2 \ \mbox{and} \ H(w_3^{\prime})=n^2+n-1.$$ Thus,  by the theorem \ref{autointer}, $$i(\gamma_1^{\prime};\gamma_1^{\prime})=i(\gamma_3^{\prime};\gamma_3^{\prime})=nL-n^2.$$ This gives the proof.

\eproof

From these corollaries we deduce:
\bcor

Let $\gamma$ be a closed geodesic of $P$, $L(\gamma)$ his combinatorial length and $k$ an integer. 
\begin{center}
If $i(\gamma;\gamma)=k$, then $2\sqrt{k}\leq L(\gamma)\leq 2k+2$.
\end{center}
\ecor

\brmq
This corollary tells us an important thing on the pair of pants; if we fix the number of self-intersections $k$, there exists an integer $L_0$ such that for every closed geodesic $\gamma$ of length $L>L_0$, his number of self-intersections is bigger than $k$. The result of M.Mirzakhani (\cite{Mir}) shows that for any hyperbolic compact surface different from the pair of pants, this result is false. 

\ermq

\section{Proof of theorem \ref{thm5}}
Now we are interested on the following set:
$$\A_{\epsilon}(L)=\{\gamma\in \mathcal{G}^c | L(\gamma)=L ,\ i(\gamma;\gamma)\geq \ds (\ds\frac{1}{4}-\epsilon)L^2 \}$$ for $\epsilon>0$.

Let $\gamma$ be a closed geodesic of $P$ and let $w=w(\gamma)=s_1^{i_1}r_1^{j_1}\cdots s_n^{i_n}r_n^{j_n}$ be the word associated to $\gamma$. By using the theorem \ref{thm1}, we can deduce the following result:
\begin{center}
{\it If the closed geodesic $\gamma\in\A_{\epsilon}(L)$ then $n\geq\ds\frac{L}{2}-\ds\sqrt{\epsilon}L$.}
\end{center}
 Set $B_\epsilon(L)=\{\gamma\in \mathcal{G} ; \ L(\gamma)=L ;\  L-2n\leq 2\sqrt{\epsilon}L\}$,
 we have the following inclusion: 
  $$A_\epsilon(L)\subset  B_\epsilon(L).$$   
 
 Let $\gamma$ be a closed geodesic of $P$ and $w(\gamma)=s_1^{i_1}r_1^{j_1}\cdots s_n^{i_n}r_n^{j_n}$ the word associated to $\gamma$.
There exists a permutation of the powers of the letters of the word $w$ and a positive integer $N$ such that $w(\gamma)$ can be write on the following form: $w=s_1^{i_1}r_1^{j_1}\cdots  s_N^{i_N}r_N^{j_N}s_{N+1}r_{N+1}\cdots s_nr_n$ with $i_k>1$ or $j_k>1$ for any $1\leq  k\leq N$.
Set $L_N=\ds\sum_{k=1}^{N}i_k+j_k$, then $L=L_N+2(n-N)$ and $L_N-2N\leq 2\sqrt{\epsilon}L.$ We have: $$\#B_\epsilon(L)=\#\{w=s_1^{i_1}r_1^{j_1}\cdots s_N^{i_N}r_N^{j_N}s_{N+1}r_{N+1}\cdots s_nr_n; \ L_N-2N\leq 2\sqrt{\epsilon}L\}.$$

For any $1\leq k\leq N$, $i_k>1$ or $j_k>1$, then $L_N\geq 3N$ and so $L_N\leq 6\sqrt{\epsilon}L$.
Thus  $$\#B_\epsilon(L)\leq \#\{w=s_1^{i_1}r_1^{j_1}\cdots s_N^{i_N}r_N^{j_N}s_{N+1}r_{N+1}\cdots s_nr_n; \ L_N\leq 6\sqrt{\epsilon}L\}.$$
A straightforward calculation gives us:
 $$\#B_\epsilon(L)\leq \ds\frac{1}{L}\sum_{L_N=0}^{E(6\sqrt{\epsilon}L)+1}8\ds(3^{L_N-2}-2^{L_N-2})(2^{L-L_N-2})$$

and we deduce the following inequality:
$$\#\A_{\epsilon}(L)\leq \ds\frac{2^L}{L}\big[2^2(\ds\frac{3}{2})^{\ds 6\sqrt{\epsilon}L}-3\sqrt{\epsilon}L\big].$$
Recall that $\#\G^{c}(L)=\ds\frac{8\times 3^{L-2}}{L}$. When $\ds\epsilon<\frac{1}{36}$, this inequality implies that: $$\ds\lim\limits_{L\to +\infty}\frac{\#\A_{\ds\epsilon}(L)}{\#\G^{c}(L)}=0.$$
 
\bpro  Let $N$ be the integer defined above and $p$ a positive integer.
Let $\gamma$ be the closed geodesic of $P$ associated to the word \\ $w=s_1^{i_1}r_1^{j_1}\cdots  s_N^{i_N}r_N^{j_N}s_{N+1}r_{N+1}\cdots s_{N+p}r_{N+p}(a\bar{b})^{n-N-p}$ with $i_k>1$ or $j_k>1$ for any $1\leq  k\leq N$.

If $L-2n+2p\leq 2\epsilon L$ then the closed geodesic $\gamma$ lies in $\A_{\epsilon}(L)$.

\epro

\bproof To prove this proposition we will use the following result: $$i(\gamma;\gamma)=L-2n+H(w)+2\ds\sum_{k=1}^{n}\sum_{l=k+1}^{n}\min(i_k;i_l)+\min(j_k;j_l)-2.$$
 We have:
$$\ds\sum_{k=1}^{n}\sum_{l=k+1}^{n}\min(i_k;i_l)+\min(j_k;j_l)-2= \ds\sum_{k=1}^{N}\sum_{l=k+1}^{N}\min(i_k;i_l)+\min(j_k;j_l)-2.$$
The definition of $N$ implies that $\min(i_k;i_l)+\min(j_k;j_l)\geq 3$ for any integers $k$ and $l$ between $1$ and $N$. This implies that: $$i(\gamma;\gamma)\geq H(w)+N^2+L_N-3N,$$ where $L_N=\ds\sum_{k=1}^{N}i_k+j_k$.
The conditions on  $p$ and $N$ give: $$H(w)\geq (n-p-N)^2+2nN.$$
As $n-p\geq \ds\frac{L}{2}-\epsilon L$, we have $i(\gamma;\gamma)\geq (\ds\frac{1}{4}-\epsilon)L^2.$
 This achieve the proof.

\eproof
This proposition shows that for every word $w=w^{\prime}a\bar{b}\cdots a\bar{b}$ such that $w^{\prime}$ is a finite word of length $L(\gamma(w^{\prime}))\leq 6 \ds\epsilon L$, the closed geodesic $\gamma$ associated to this word, lies in $\A_\epsilon(L)$. This implies that:
$$ \#\{w=w^{\prime}(a\bar{b})^{n-p-N}; \ L(\gamma(w^{\prime}))\leq 6\epsilon L\}\leq \#\A_\epsilon(L)$$ and so
 $$\#\A_\epsilon(L)\geq \ds\frac{2}{3L}[3^{6\epsilon L}-1].$$

This gives the theorem \ref{thm5}.


\begin{paragraph}{Acknowledgement:} The authors would like to thank the project NLAGA for their support.
\bibliographystyle{plain}
\end{paragraph}
\bibliography{self_intersection_on_a_pair_of_pants}

\textbf{Departement of mathematics, Universit\'e  Cheikh Anta Diop, Dakar, Senegal.}\\
\textit{email: elhadjiabdouaziz2.diop@ucad.edu.sn}

\noindent \textbf{Departement of mathematics, Universit\'e  Cheikh Anta Diop, Dakar, Senegal.}\\
\textit{email: masseye.gaye@ucad.edu.sn}

\end{document}